\documentclass[12pt]{amsart}

\usepackage{verbatim}

\newtheorem{theorem}{Theorem}[section]

\newtheorem{corollary}[theorem]{Corollary}

\newtheorem{lemma}[theorem]{Lemma}

\newtheorem{proposition}[theorem]{Proposition}

\theoremstyle{definition}

\newtheorem{definition}[theorem]{Definition}

\newtheorem{remark}[theorem]{Remark}

\theoremstyle{parrafo}

\newtheorem{parrafo}[theorem]{{\!}}

\numberwithin{equation}{theorem}

\newcommand{\calo}{{\mathcal {O}}}

\newcommand{\Sing}{\mbox{Sing\ }}

\title[]{Rees algebras on smooth
schemes: integral closure and higher differential operator.}

\author{Orlando Villamayor U.}

\address{Dpto. Matem\'aticas, Facultad de Ciencias, Universidad
Aut\'onoma de Madrid, Canto Blanco 28049 Madrid, Spain.}
\email{villamayor@uam.es}
\thanks{.}


 \date{May 2006}
\dedicatory{} \commby{}

\begin{document}
\maketitle
\begin{abstract}
Let $V$ be a smooth scheme over a field $k$, and let $\{I_n, n\geq
0\}$ be a filtration of sheaves of ideals in $\calo_V$, such that
$I_0=\calo_V$, and $I_s\cdot I_t\subset I_{s+t}$. In such case
$\bigoplus I_n$ is called a Rees algebra.

A Rees algebra is said to be a Diff-algebra if, for any
two integers $N>n$ and any differential operator $D$ of order $n$,
$D(I_N)\subset I_{N-n}$. Any Rees algebra extends to a
smallest Diff-algebra.

There are two ways to define extensions of Rees algebras, and both are of interest in singularity theory. One is that defined by taking
integral closures (in which a Rees algebra is included in its integral closure), and another extension is that defined, as above, in which the algebra is extended to a Diff-algebra.

Surprisingly enough, both forms of extension are compatible in a natural way. Namely, there is a compatibility of higher differential operators with integral closure which we explore here under the assumption that $V$ is smooth over a perfect field.


\end{abstract}

\tableofcontents






\label{introduction}

\section{Introduction}

A smooth ring $R$, of finite type over a field $k$, has a locally
free sheaf of $k$-linear differential of order $n$, say
$Diff^n(R)$, for each index $n\geq 0$.

An $R$ sub-algebra of a polynomial ring over $R$, say $R[W]$, can be
expressed as $\mathcal{G}=\bigoplus_{n\geq 0}I_nW^n$, where $I_0=R$, and each
$I_n$ is an ideal in $R$. 
$\mathcal{G}$ is called the Rees algebra
of the filtration ${\{I_n\}}_{n\geq 0}$. For example a Rees ring
of an ideal $I$ is of this type; in which case $I_n=I^n$. Notice that the 
integral closure of a Rees ring is also a Rees algebra (not
necessarily a Rees ring of an ideal).

Taking integral closure of Rees algebras included in $R[W]$ can be thought
of as an operator, say $\mathcal{G}\subset
\widetilde{\mathcal{G}}$

The study of invariants for embedded singularities has motivated another kind of
extension, linked to differential operators. In fact, from the
point of view of singularities it is interesting to consider Rees
algebras $\mathcal{G}=\bigoplus_{n\geq 0}I_nW^n$ with an
additional compatibility with differential operators. Namely with the 
property that for any operator $D\in Diff^n(R)$, and
any index $N>n$, $D(I_N)\subset I_{N-n}$.

Rees algebras with this property are called here Diff-algebras, and we show that for any $\mathcal{G}=\bigoplus_{n\geq
0}I_nW^n$, there is a smallest extension to one with this
property, say $G(\mathcal{G})$. This defines a second operator,
say $\mathcal{G}\subset G(\mathcal{G})$. 
The objective of this paper is to study the interplay
between both operators. In Main Theorem \ref{TH515} it is shown that if two
Rees algebras have the same integral closure, then their
$G$-extensions also have the same integral closure. This shows a
curious relation of differential operators with integral closures.

Sections 2 and 3 are devoted to extensions of Rees algebras to those which are compatible with differential operators. The main result in the latter section is Theorem \ref{3th1} that characterizes the extension $\mathcal{G}\subset G(\mathcal{G})$.

 In Sections  4 and 5 we study some natural relation of the subjects discussed in the first two sections with singularities.
 The main results in Section 5 is Theorem \ref{th54}. Rees algebras compatible with differential operators are studied here only for algebras over smooth schemes. This Theorem studies this compatibility when a morphism is defined between two smooth schemes.
  
 
 Our main results relating integral closure of Rees algebras and differential operators are finally addressed in Section 6. 

This is a paper on commutative algebra, and the outcome has lead to
results in singularities theory in positive characteristic, most particularly in the elimination theory developed with this purpose in (\cite{VV2}).

The study of differential operatores on smooth schemes, and their relation with with integral closure has also been studied,  independently, in a very interesting paper of Kawanoue. The presentation in \cite{kaw} already includes some results treated here.

Extensions of algebras by differential operators appeared  in   \cite{You} , and more recently in work of
 J. Wlodarczyk , and J. Koll\'{a}r  (\cite{kollar}, and \cite{WL}).  Also in,   \cite{Hironaka77},\cite{Hironaka03},\cite{Hironaka05}
where Hironaka studies the interplay of differential operators with the theory of infinitely close points.
Our presentation here, although related with these last three papers, makes no use of monoidal transformations, or of
infinitely close points. 



I profited from discussions with Vincent Cossart, Marco Farinati,
Monique Lejeune, Luis Narv\'{a}ez, and Santiago Zarzuela.

\section{ Graded rings and Diff-algebras}

\begin{parrafo}\label{IIpar1}  {\em Fix a noetherian ring $B$ and a sequence of ideals $\{I_k\}$, $k\geq 0$,
which fulfill the conditions:

{\bf 1)} $I_0=B$, and

{\bf 2)} $I_k\cdot I_s \subset I_{k+s}$.

This defines a $B$-algebra which is a graded subring $\mathcal{G}=\bigoplus_{k\geq 1}
I_k W^k$ of the polynomial ring $ B[W]$.  We say that
$\mathcal{G}$ is a {\em Rees algebra}  if this subring is a
(noetherian) finitely generated $B$-algebra.

\vskip 0.5cm

In what follows we define a Rees algebra, say
$\bigoplus_{n\geq 0} I_nW^n$ in $B[W]$, by fixing a set of
generators, say $$\mathcal{F}=\{g_{n_i}W^{n_i}/ g_{n_i}\in B, n_i>0, 1\leq i \leq
m\}.$$ 
So if $f \in I_n$, then $$f=F_n(g_{n_1},\dots , g_{n_m}),$$
where $F_n(Y_1,\dots , Y_m)$ is a weighted homogeneous polynomial
in $m$ variables, each $Y_j$ considered with weight $n_j$.}
\end{parrafo}

\begin{remark}\label{IIrk1}

1) Examples of Rees algebras are  Rees rings of  ideals, say 
$I\subset B$. In this case $I_k=I^k$ for each $k\geq 1$. These are algebras generated by  (homogenoeus) elements of degree one (i.e., generated by $\mathcal{F}$ with all $n_i=1$).

2) When $\bigoplus I_k W^k \subset (\subset B[W])$ is a
Rees algebra,  a new Rees algebra $\bigoplus
I'_k W^k $ is defined by setting $$ I'_k=\sum_{r\geq k} I_r.$$
Note that $I'_k \supset I'_{k+1}$.

If $\bigoplus I_kW^k  $ is generated by $\mathcal{F}=\{
g_{n_i}W^{n_i}, n_i>0, 1\leq i\leq m \}$. Namely, if: $$ \bigoplus
I_kW^k =R[\{g_{n_i}W^{n_i}\}_{g_{n_i}W^{n_i}\in
\mathcal{F}}],$$ then we claim that  

A) $\bigoplus I'_k W^k $ is generated by
the finite set $\mathcal{F}'=\{ g_{n_i}W^{n'_i}, 1\leq i\leq m, 1 \leq n'_i\leq
n_i \}$, and

B) $\bigoplus I_k
W^k \subset \bigoplus I'_k W^k $ is a finite extension.

To prove the first claim we can use the fact that an element, say $h_N$, is homogeneous of degree $N$ in the $B$-subalgebra generated by 
 $\mathcal{F}'$ if it is a $B$-linear combination  of monomials of the form $$h_1^{a_1}\cdot h_2^{a_2} \cdots h_s^{a_s}$$
 where $h_iW^{n_i}\in \mathcal{F}'$, and $\sum a_i\cdot n_i=N$. Suppose that $a_1\neq 0$, and express $$h_1^{a_1}\cdot h_2^{a_2} \cdots h_s^{a_s}= h_1\cdot h_1^{a_1-1}\cdot h_2^{a_2} \cdots h_s^{a_s}$$ where now the first factor  $h_1$ is endowed with degree $n_1-1$. This ensures that $h_1^{a_1}\cdot h_2^{a_2} \cdots h_s^{a_s}$ also appears, in the Rees algebra, as an homogeneous element in degree 
 $N-1$ (as an element in $I_{N-1}$). This already proves the first claim.

To prove B), it suffices to check that given $g\in I_{k}$,
then $g W^{k-1}$ is integral over $\bigoplus I_kW^k $.
Note that $$g\in I_{k}\Rightarrow g^{k-1} \in
I_{k(k-1)}\Rightarrow g^{k} \in I_{k(k-1)},$$ so $g W^{k-1}$
fulfills monic polynomial the equation $Z^{k}-(g^{k}W^{k(k-1)})=0$.

Here we always assume that $B$ is an excellent ring, so that the integral closure is also finitely generated over $B$. So B) shows that up to integral closure we may assume that a Rees algebra has the
additional condition:  $$I_k\supset I_{k+1}.$$

\end{remark}

\begin{parrafo}\label{IIpar3}
{\em 

For a fixed positive integer $N$, $B[W^N]\subset B[W]$ is a finite
extension of graded rings. Furthermore, $\bigoplus_{k\geq 0}
I_{kN}W^{kN}$ is a Rees algebra in $B[W^N]$, and the ring extension $\bigoplus_{k\geq
0} I_{kN}W^{kN}\subset \bigoplus_{n\geq 0} I_nW^n$ is also finite.

If the Rees algebra is generated by $\mathcal{F}=\{g_{n_i}W^{n_i}/ n_i>0, 1\leq i \leq
m\}$, and if $N$ is a common multiple of all integers $n_i$,
$1\leq i\leq m$, then $$\bigoplus_{k\geq 0} I_{N}^kW^{kN}\subset
\bigoplus_{n\geq 0} I_nW^n$$ is an integral extension. Here the
left hand side is the Rees ring of $I_N$ (in $B[W^N]$). So Rees
algebras are, at least in this sense, finite extensions of Rees rings.

If a Rees algebra $\bigoplus_{n\geq 0} I_nW^n$ in $B[W]$ is the
Rees ring of $I_1$, then the integral closure in $B[W]$ is
$\bigoplus_{n\geq 0} \overline{I}_nW^n$, where each
$\overline{I}_n$ is the integral closure of the ideal $I_n$. This
is a Rees algebra, and not necessarily the Rees ring of the ideal
$\overline{I}_1$.}
\end{parrafo}
\begin{parrafo}\label{IIpar31}
{\em Let $B$ be a normal excellent ring, and let
$$Spec(B)\stackrel{\pi}{\longleftarrow} X$$
be a proper birational morphism, then ${I}\subset
\pi_{*}(I\calo_X) \subset \overline{I},$ where $\overline{I}$
denotes the integral closure of $I$ in $B$. Moreover, if $\pi$ is
the normalization of the blow-up at $I$, then $I\calo_X$ is an
invertible sheaf of ideals, and
$$\overline{I}=\pi_{*}(I\calo_X).$$

Assume that the normal ring $B$ is of finite type over a field
$k$. If $B$ is a one dimensional normal domain, any ideal is
invertible and integrally closed. We add the following well known
result for self-containment ( se \cite{Hironaka77}, p.54 or \cite{LT} p. 100).

\begin{lemma}\label{lemclauent} Let $I$, $J$ be two ideals in a normal domain $B$,
which is finitely generated over a field $k$. Then
$\overline{I}=\overline{J}$ if and only if $I\calo_W=J\calo_W$,
for any morphism of $k$-schemes $W \to Spec(B)$, with $W$ of
dimension one, regular and of finite type over $k$.

\end{lemma}

\proof Let $x\in W$ be a closed point that maps to, say $y\in Spec(B)$. Then $\calo_{W,x}$ is a valuation
ring that dominates $\calo_{Spec(B), y}$. So if
$\overline{I}=\overline{J}$, then $I\calo_W=J\calo_W$. In fact,
for any morphism $B\to A$, where $A$ is a valuation ring,
$IA=\overline{I}A$.

Assume that this condition holds for any morphism from a regular
one dimensional scheme $W$. We claim now that both ideals have the
same integral closure in $B$.

Let $Spec(B)\stackrel{\pi}{\longleftarrow} X$ be the normalized
blow up at $I$, and let $\{H_1,\dots, H_s\}$ be the irreducible
components of the closed set defined by the invertible sheaf of
ideals $I\calo_X$. Here each $H_i$ is an irreducible hypersurfaces
in $X$. Let $h_i\in X$ denote the generic point of $H_i$. There
are positive integers $a_i$, so that $I\calo_X$ can be
characterized as the sheaf of functions vanishing along $H_i$ with
order at least $a_i$ (i.e., with order at least $a_i$ at the
valuation rings $\calo_{X,h_i}$).

{\em Claim: }The sheaf of ideals $J\calo_X$ also has order $a_i$
at $\calo_{X,h_i}$.

If the claim holds, $J\calo_X\subset I \calo_X$, and
$$J \subset \pi_{*}(J\calo_{X}) \subset \pi_{*}(I\calo_{X})=
\overline{I}.$$ In particular $\overline{J}\subset \overline{\overline{I}}=\overline{I}$.
A similar argument leads to the other inclusion.

In order to prove the claim we choose a closed point $x\in H_i$ so
that:

1) $\calo_{X,x}$ is regular,

2) $x\in H_i-\cap_{j\neq i} H_j$,

3) $H_i$ is regular at $x$, and

4) $J\calo_{X,x}$ is a $p$-primary ideal, for $p=I(H_i)_x$.

Since any sheaf of ideals has only finitely many $p$-primary
components, such choice of $x$ is possible.

Let $\{x_1,\dots x_{d-1},x_d\}$ be a regular system of parameters
for $\calo_{X,x}$ such that $p=I(H_i)_x=x_d\calo_{X,x}$, and let
$W$ be the closure of the irreducible curve defined locally by
$<x_1,\dots ,x_{d-1}>$. So $W$ is one dimensional, and regular
locally at $x$. We may assume that $W$ is regular after applying
quadratic transformations which do not affect the local ring
$\calo_{W,x}$. By construction $I\calo_{W,x}$ has order $a_i$, by
hypothesis the same holds for $J\calo_{W,x}$. This proves the
claim.

 }
\end{parrafo}

\begin{parrafo}\label{IIpar2}
{\em Let $B=S[X]$ be a polynomial ring, and let $Tay: B \to B[U]$
be the $S$-algebra homomorphism defined by setting $Tay(X)=X+U$.
For any $f(X)\in B$ set
$$Tay(f(X))=\sum_{\alpha \geq 0}
\Delta^{\alpha}(f(X))U^{\alpha}.$$
This defines, for each $\alpha$, $\Delta^{\alpha}:S[X]\to S[X]$, which is an $S$-differential operators
($S$ linear) on $B=S[X]$. Furthermore, for any positive integer $N$, the set
$\{\Delta^{\alpha}, 0\leq \alpha \leq N\}$ is a basis of the
$B$-module of $S$-differential operators on $B$, of order $\leq
N$.}

\end{parrafo}

\begin{definition}\label{IIdef1}
Let $B=S[X]$ be a polynomial ring over a noetherian ring $S$. A
Rees algebra $$\bigoplus I_k\cdot W^k\subset B[W]$$ is said to be a Diff-algebra, relative to $S$, when:

i) $I_k\supset I_{k+1}$ for all $k\geq 0$.

ii) For all $n>0$ and $f\in I_n$, and for any index $0\leq j \leq
n$ and any $S$-differential operator of order $\leq j$, say $D_j$:
$$D_j(f)\in I_{n-j}.$$
\end{definition}

\begin{remark}\label{IIrk2} Let $Diff^N_S(B)$ denote the module of
$S$-differential operators of order at most $N$. Then
$$Diff^N_S(B)\subset Diff^{N+1}_S(B)\subset \dots$$ For this
reason it is natural to require condition (i) in our previous
definition. Note also that \ref{IIpar2} asserts that (ii) can be
reformulated as:

\bigskip

ii') For any $n>0$ and $f\in I_n$, and for any index $0\leq \alpha
\leq n$: $$\Delta^{\alpha}(f)\in I_{n-\alpha},$$

In fact, (i) and (ii) are equivalent to (i) and (ii'):
\end{remark}
\begin{theorem}\label{IIth} Fix $B=S[X]$ as before, and a finite set $\mathcal{F}=\{
g_{n_i}W^{n_i},  n_i>0 , 1\leq i\leq m \}$, with the following
properties:

a) For any $1\leq i\leq m$, and any $n'_i$, $0< n'_i\leq n_i$:
$$g_{n_i}W^{n'_i}\in \mathcal{F}.$$

b) For any $1\leq i\leq m$, and for any index $0\leq \alpha <
n_i$: $$\Delta^{\alpha}(g_n)W^{n_i-\alpha}\in \mathcal{F}.$$

Then the $B$ subalgebra of $B[W]$, generated by $\mathcal{F}$ over
the ring $B$, is a  Diff-algebra relative to $S$.
\end{theorem}
\proof

Condition (i) in Def \ref{IIdef1} holds by \ref{IIrk1}, 2).

Fix a positive integer $N$, and let $I_NW^N$ be the homogeneous component of degree $N$ of the $B$
subalgebra generated by $\mathcal{F}$. We prove that for any $h\in
I_N$, and any $0\leq \alpha \leq N$, $\Delta^{\alpha}(h)\in
I_{N-\alpha}$.

The ideal $I_N \subset B$ is generated by all elements of the form
\begin{equation}\label{IIeq1}
H_N=g_{n_{i_1}}\cdot g_{n_{i_2}}\cdots g_{n_{i_p}} \hskip 1cm
n_{i_1}+n_{i_2}+\cdots n_{i_p}=N,
\end{equation}
with the $g_{n_{i_i}}W^{n_{i_i}}\in \mathcal{F}$ not necessarily
different.

Since the operators $\Delta^{\alpha}$ are linear, it suffices to
prove that $\Delta^{\alpha}(a\cdot H_N)\in I_{N-\alpha}$, for
$a\in B$, $H_N$ as in \ref{IIeq1}, and $0\leq \alpha \leq N$. We
proceed in two steps, by proving:

1)  $\Delta^{\alpha}(H_N)\in I_{N-\alpha}$.

2)  $\Delta^{\alpha}(a\cdot H_N)\in I_{N-\alpha}$.

We first prove 1). Set $Tay: B=S[X] \to B[U]$, as in \ref{IIpar2}.
Consider, for any element $g_{n_{i_l}}W^{n_{i_l}}\in \mathcal{F}$,
$$Tay(g_{n_{i_l}})=\sum_{\beta \geq 0}
\Delta^{\beta}(g_{n_{i_l}})U^{\beta}\in B[U].$$

Hypothesis (b) states that for each index $0\leq \beta < n_{i_l}$,
$\Delta^{\beta}(g_{n_{i_l}})W^{n_{i_l}-\beta}\in \mathcal{F}$.

On the one hand $$Tay(H_N)=\sum_{\alpha \geq 0}
\Delta^{\alpha}(H_N)U^{\alpha},$$ and, on the other hand
$$Tay(H_N)=Tay(g_{n_{i_1}})\cdot Tay(g_{n_{i_2}})\cdots
Tay(g_{n_{i_p}})$$ in $B[U]$. This shows that for a fixed $\alpha$
($0\leq \alpha \leq N$), $\Delta^{\alpha}(H_N)$ is a sum of
elements of the form: $$\Delta^{\beta_1}(g_{n_{i_1}})\cdot
\Delta^{\beta_2}(g_{n_{i_2}})\cdots
\Delta^{\beta_p}(g_{n_{i_p}}),\hskip 1cm \sum_{1\leq s\leq p}
\beta_s=\alpha.$$

So it suffices to show that each of these summands is in
$I_{N-\alpha}$.

Note here that $$\sum_{1\leq s\leq p}
(n_{i_s}-\beta_s)=N-\alpha,$$ and that some of the integers
$n_{i_s}-\beta_s$ might be zero or negative. Set $$G=\{ r, 1\leq r
\leq p, \mbox{ and } n_{i_r}-\beta_r > 0\}.$$ So
$$N-\alpha=\sum_{1\leq s\leq p} (n_{i_s}-\beta_s)\leq \sum_{r\in
G} (n_{i_r}-\beta_r)=M.$$ Hypothesis (b) ensures that
$\Delta^{\beta_r}(g_{n_{i_r}})\in I_{n_{i_r}-\beta_r}$ for every
index $r\in G$, in particular:$$\Delta^{\beta_1}(g_{n_{i_1}})\cdot
\Delta^{\beta_2}(g_{n_{i_2}})\cdots
\Delta^{\beta_p}(g_{n_{i_p}})\in I_M.$$ Finally, since $M \geq
N-\alpha,$ $I_M\subset I_{N-\alpha}$, and this proves Case 1).

For Case 2), fix $0\leq \alpha \leq N$. We claim that
$\Delta^{\alpha}(a\cdot H_N)\in I_{N-\alpha}$, for $a\in B$ and
$H_N$ as in \ref{IIeq1}. At the ring $B[U]$, $$Tay(a\cdot
H_N)=\sum_{\alpha \geq 0} \Delta^{\alpha}(a \cdot
H_N)U^{\alpha},$$ and, on the other hand $$Tay(a \cdot
H_N)=Tay(a)\cdot Tay(H_N).$$ This shows that
$\Delta^{\alpha}(a\cdot H_N)$ is a sum of terms of the form $
\Delta^{\alpha_1}(a)\cdot \Delta^{\alpha_2}(H_N)$,$ \alpha_i\geq
0$, and $\alpha_1+\alpha_2=\alpha$. In particular $\alpha_2 \leq
\alpha$; and by Case 1), $\Delta^{\alpha_2}(H_N)\in
I_{N-\alpha_2}$. On the other hand $N-\alpha_2\geq N-\alpha$, so
$\Delta^{\alpha_2}(H_N)\in I_{N-\alpha}$, and hence
$\Delta^{\alpha}(a\cdot H_N)\in I_{N-\alpha}$.
\endproof
\begin{corollary}\label{2corth} The Rees algebra in $B[W]$, generated over
$B$ by $$\mathcal{F}=\{ g_{n_i}W^{n_i},  n_i>0 , 1\leq i\leq m
\},$$ extends to a smallest Diff-algebra, which is generated by
the finite set
$$\mathcal{F'}=\{\Delta^{\alpha}(g_n)W^{n_i-\alpha}/
g_{n_i}W^{n'_i}\in \mathcal{F}, \mbox{ and } 0\leq \alpha <
n_i\}.$$

\end{corollary}

\begin{remark} \label{rky} (Not used in what follows) Theorem \ref{IIth} shows how to extend any Rees algebra to a Diff-algebra, say 
$\bigoplus I_kW^k\subset B[W]$ so that the conditions of Definition \ref{IIdef1} holds; namely that 
for any $S$-differential operator of order  $ j(\leq n)$, say $D_j$: $D_j(I_n)\in I_{n-j}.$

A similar argument can be used to extend Rees algebras to algebras, say $\bigoplus I_kW^k\subset B[W]$ again,  with the condition :
\begin{equation}\label{eq23}
D_j(I_n)\in I_{n}
\end{equation}
for any positive $n$, and any differential operator of order $j$, with no condition on $j$. It is easy to check that ideals $I_n$ with this property are those generated by elements in $S$.

Consider, as in Theorem \ref{IIth}, a finite set $\mathcal{F}=\{
g_{n_i}W^{n_i},  n_i>0 , 1\leq i\leq m \}$, with the following
properties:

a) For any $1\leq i\leq m$, and any $n'_i$, $0< n'_i\leq n_i$:
$g_{n_i}W^{n'_i}\in \mathcal{F}.$

b) For any $1\leq i\leq m$, and for any index $0\leq \alpha $: $\Delta^{\alpha}(g_{n_i})W^{n_i}\in \mathcal{F}.$

We claim now that the $B$ subalgebra of $B[W]$, generated by $\mathcal{F}$ over
the ring $B$, fulfills (\ref{eq23}). 
Note here that each $g_{n_i}$ is polynomial in $X$, so $\Delta^{\alpha}(g_n)=0$ for $\alpha$ big enough, so  $\mathcal{F}$ is in fact finite.

In order to prove the claim it suffices to show that   $\Delta^{\alpha}(a\cdot H_N)\in I_{N}$, for
$a\in B$, $H_N$ as in \ref{IIeq1}. As in the previous Theorem we
proceed in two steps, but proving now that:

1)  $\Delta^{\alpha}(H_N)\in I_{N}$.

2)  $\Delta^{\alpha}(a\cdot H_N)\in I_{N}$.

$\Delta^{\alpha}(H_N)$ is a sum of
elements of the form: $\Delta^{\beta_1}(g_{n_{i_1}})\cdot
\Delta^{\beta_2}(g_{n_{i_2}})\cdots
\Delta^{\beta_p}(g_{n_{i_p}}),\hskip 0.1 cm \sum_{1\leq s\leq p}
\beta_s=\alpha.$ So, to prove 1),  it suffices to show that each of these products is in
$I_{N}$. This follows from (\ref{IIeq1}) and the assumption on  $\mathcal{F}$. The proof for 2) is similar.

\end{remark}

\begin{remark} \label{rky1} (Not used in what follows.) The proof of Theorem \ref{IIth}, and also that of Remark \ref{rky}, rely strongly 
on the fact that $Tay: B \to B[U]$, defined on the polynomial ring  $B=S[X]$ by setting $Tay(X)=X+U$, is an $S$-algebra homomorphism. In fact the proof of the Theorem reduce  to showing that $\Delta^{\alpha}(H_N)\in I_{N-\alpha}$ (that $\Delta^{\alpha}(H_N)\in I_{N-\alpha}$ in the case of Remark \ref{rky}), where $H_N=g_{n_{i_1}}\cdot g_{n_{i_2}}\cdots g_{n_{i_p}}$ is a product of elements in a finite set of generators $\mathcal{F}$.

An interesting alternative $S$-algebra homomorphism is $$Tay_X:B \to B[U],$$
defined  by setting $Tay_X(X)=X+XU$. In this case 
$$ Tay_X(F(X)=\sum_{\alpha \geq 0}
X^{\alpha}\Delta^{\alpha}(f(X))U^{\alpha}.$$

If a finite set $\mathcal{F}=\{
g_{n_i}W^{n_i},  n_i>0 , 1\leq i\leq m \}$ is such that:

a) For any $1\leq i\leq m$, and any $n'_i$, $0< n'_i\leq n_i$:
$g_{n_i}W^{n'_i}\in \mathcal{F}.$

b) For any $1\leq i\leq m$, and for any index $0\leq \alpha $: $X^{\alpha}\Delta^{\alpha}(g_{n_i})W^{n_i}\in \mathcal{F}.$

As each $g_{n_i}$ is polynomial in $X$,  $X^{\alpha}\Delta^{\alpha}(g_n)=0$ for $\alpha$ big enough, so  $\mathcal{F}$ is in fact finite.

The same argument used above show that for these algebras:
$$X^{\alpha}\Delta^{\alpha}(I_n)\subset I_n.$$
Rees algebras with this property are considered in toric geometry. They are also characterized by the fact that if $f(X)=\sum s_r X^r(\in S[X])$ is in $I_n$, then each  $s_r X^r\in I_n$.

\end{remark} 




\section{ Diff-algebras over smooth schemes.}

\begin{parrafo}\label{par3.0}

{\rm A sequence of coherent ideals on a scheme $Z$, say
$\{I_n\}_{n\in \mathbb{N}}$, such that $I_0=\mathcal{O}_Z$, and
$I_k\cdot I_s\subset I_{k+s}$, defines a graded sheaf of algebras
$\bigoplus_{n\geq 0} I_nW^n \subset  \mathcal{O}_{Z}[W]$.

We say that this algebra is a Rees algebra if there is an open
covering of $Z$ by affine sets $\{U_i\}$, so that each restriction
$$\bigoplus_n I_n(U_i) W^n \subset \calo_Z(U_i)[W]$$ is a finitely generated
$\mathcal{O}_Z(U_i)$-algebra.

In what follows $Z$ will denote a smooth scheme over a perfect field $k$,
and $Diff^r_k(Z)$, or simply $Diff^r_k$, the locally free sheaf of
$k$-linear differential operators of order at most $r$. 
}

\end{parrafo}

\begin{definition} \label{3def1} We say that a Rees algebra defined by $\{I_n\}_{n\in
\mathbb{N}}$ is a Diff-algebra relative to the field $k$, if:

i) $I_n\supset I_{n+1}$.

ii) There is open covering of $Z$ by affine open sets $\{U_i\}$,
and for any $D\in Diff^{(r)}(U_i)$, and any $h\in I_n(U_i)$, then
 $D(h)\in I_{n-r}(U_i)$, provided $n\geq r$.

\end{definition}

Due to the local nature of the definition, we reformulate it in
terms of smooth $k$-algebras.

\begin{definition}
In what follows $R$ will denote a smooth algebra over a perfect field, or
a localization of such algebra at a closed point ( a regular local
ring). A Rees algebra is defined by a sequences of ideals
$\{I_k\}_{k\in \mathbb{N}}$ such that:

1) $I_0=R$, and $I_k\cdot I_s \subset I_{k+s}$.

2) $\bigoplus I_k W^k$ is a finitely generated $R$-algebra.

We shall say that the Rees algebra is a  Diff-algebra relative to
$k$, if

3) $I_n\supset I_{n+1}$, and

4) given $D\in Diff_k^{(r)}(R)$, then $D(I_n)\subset I_{n-r}$.

\end{definition}

\bigskip

We now show that any Rees algebra extends to a smallest
Diff-algebra  (i.e., included in any other Diff-algebra
containing it).

\begin{theorem}\label{3th1}
Fix a smooth scheme $Z$ over a perfect field $k$. Assume that $\mathcal{G}=\bigoplus I_k W^k$ is a Rees algebra
over $Z$. Then there is a natural and smallest
extension of it, say $\mathcal{G}\subset G(\mathcal{G})$,  where $G(\mathcal{G})$ is a Diff-algebra relative to the field $k$.
\end{theorem}
\proof The problem is local, so we will assume that $R$ is the
local ring at a closed point, and show that a finitely generated
subalgebra of $R[W]$ extends, by successive applications of
differential operators, to a finitely generated algebra.

We will argue in steps. Assume that the local ring $R$ is of
dimension 1, and let $x$ denote a parameter. Set $Tay: \hat{R} \to
\hat{R}[[U]]$ the k-algebra morphism at the completion defined by
setting $Tay(x)=x+U$. Here $\hat{R}=k'[[x]]$ is a ring of formal
power series over a finite extension $k'$ of $k$, $$Tay(f)=\sum
\Delta^r(f(x)) U^r,$$ and each $$\Delta^r:k'[[x]] \to k'[[x]]$$
maps $R$ into $R$, defining $$\Delta^r:R \to R.$$ So $Tay: \hat{R}
\to \hat{R}[[U]]$  induces by restriction $Tay:{R} \to R[[U]]$.

For any $f\in R$ set
$$Tay(f)=\sum_{r \geq 0} \Delta^{r}(f)U^{r} (\in
R[[U]]).$$

The operators $\Delta^{r}$, $r\geq 0$, are a basis of the $k$-linear
differential operators on $R$.

The same argument used in Theorem \ref{IIth} shows that if
$\bigoplus I_k\cdot W^k$ is generated by $\mathcal{F}=\{
g_{n_i}W^{n_i},  n_i>0 , 1\leq i\leq m \}$, then
$$\mathcal{F'}=\{\Delta^{r}(g_n)W^{n'_i-r}/
g_{n_i}W^{n_i}\in \mathcal{F}, \mbox{ and } 0\leq r < n'_i \leq
n_i\}$$ generates the smallest extension to a Diff-algebra.

Let now $R$ be a localization of an arbitrary smooth algebra at a
closed point, and fix a regular system of parameters $\{x_1,\dots
,x_n\}$. Define $$Tay: \hat{R} \to \hat{R}[[U_1,\dots ,U_n]]$$ as
the continuous morphisms of algebras defined by setting
$Tay(x_i)=x_i+U_i$. So for any $h\in \hat{R} $ set:
$$Tay(h)=\sum_{\alpha \in (\mathbb{N})^n}
\Delta^{\alpha}(h)U^{\alpha}.$$

This morphism defines, by restriction, $Tay: {R} \to
{R}[[U_1,\dots ,U_n]]$, and we set $$Tay(g)=\sum_{\alpha \in
(\mathbb{N})^n} \Delta^{\alpha}(g)U^{\alpha}.$$

The assumption that $k$ is perfect ensures that $\{ \Delta^{\alpha} /\alpha \in (\mathbb{N})^n, 0\leq |\alpha|\leq
n\}$ is a basis of the free $R$-module $Diff^n(R)$, and in order
to show that a Rees algebra $\bigoplus I_k\cdot W^k$ is a 
Diff-algebra, it suffices to check that given $g\in I_m$:

\begin{equation}\label{3eq1}
\Delta^{\alpha}(g)\in I_{m-|\alpha|}.
\end{equation}

Note that
$\Delta^{\alpha}\Delta^{\alpha'}=\Delta^{\alpha'}\Delta^{\alpha}.$
Define, for each index $i_0$, $1\leq i_0 \leq n$: $$ Tay_{i_0}:R
\to R[[U_{i_0}]],$$ $Tay_{i_0}(x_j)=x_j$ and
$Tay_{i_0}(x_{i_0})=x_{i_0}+U_{i_0}.$ So $$Tay_{i_0}(g)=\sum_{\alpha
\in \mathbb{N}} \Delta_{i_0}^{\alpha}(g)U^{\alpha},$$ is defined
in terms of the differential operators
$\Delta_{i_0}^{\alpha}$.
For any  $\alpha=(\alpha_1, \alpha_2,\dots , \alpha_n)
\in (\mathbb{N})^n$: $$\Delta^{\alpha}= \Delta_1^{\alpha_1} \cdots
\Delta_n^{\alpha_n},$$ is a composition of partial operators
defined above. And $\bigoplus I_kW^k$ is a  Diff-algebra if
the requirement in (\ref{3eq1}) holds for each of these partial
differential operators.

So again, the arguments in Theorem \ref{IIth} ensure that if $\bigoplus I_kW^k$ is generated by
$$\mathcal{F}=\{ g_{n_i}W^{n_i}, n_i>0 , 1\leq i\leq m \},$$ then
\begin{equation}\label{eq3422} \small
\mathcal{F'}=\{\Delta^{\alpha}(g_n)W^{n'_i-\alpha}/
g_{n_i}W^{n_i}\in \mathcal{F}, \alpha=(\alpha_1, \alpha_2,\dots ,
\alpha_n) \in (\mathbb{N})^n, \mbox{ and } 0\leq |\alpha| < n'_i
\leq n_i\}
\end{equation}
generates the smallest extension of $\bigoplus I_kW^k$ to a
Diff-algebra relative to the field $k$.

\endproof

\begin{remark}\label{rky2} (Not used in what follows.) In the previous discussion we reduce the proof of the Theorem to the case of one variable, and we make use of Theorem \ref{IIth}. There are interesting variations in the one variable case discussed  in Remark \ref{rky1}, of particular interest in the case of differentials with logarithmic poles. Such is the case when we fix an integer $s$, $1\leq s \leq n$, and consider, for each 
index $ 1\leq i_0 \leq s$, the modified function: 
$$ Tay_{x_{i_0}}:R \to R[[U_{i_0}]],$$ 
defined by $Tay_{x_{i_0}}(x_j)=x_j$ and
$Tay_{x_{i_0}}(x_{i_0})=x_{i_0}+x_{i_0}U_{i_0}.$ 

There is an natural analog of Diff-algebras with Rees algebras which are closed by differential operatores with logarithmic poles. This follows from Remark \ref{rky1}, and it is simple to extend the outcome of (\ref{eq3422})
to this context.
\end{remark}

\begin{corollary}\label{coragre} Given inclusions of
Rees algebras, say $$\mathcal{G}=\bigoplus I_n W^n \subset
\mathcal{G'}=\bigoplus I'_n W^n \subset G(\mathcal{G})=\bigoplus
I''_n W^n, $$ where $G(\mathcal{G})$ is the Diff-algebra spanned
by $\mathcal{G}$, then $G(\mathcal{G}) $ is also the Diff-algebra
spanned by $\mathcal{G}'$.

\end{corollary}
\begin{parrafo}\label{par3rela}
{\rm Fix now a smooth morphism of smooth schemes, say $Z\to Z'$.
Let $Diff^r_{Z'}(Z)$, or simply $Diff^r_{Z'}$ denote the locally
free sheaf of relative differential operators of order $r$.

We say that the Rees algebra $\bigoplus I_kW^k$
 over $Z$ (\ref{par3.0}) is a  Diff-algebra relative to $Z'$, if
conditions in Def \ref{3def1} hold, where we now require that
$D\in Diff_{Z'}^{(r)}(U_i)$ in (ii).

Since $Diff^r_{Z'}(Z)\subset Diff^r_{k}(Z)$ it follows that any
Diff-algebra relative to $k$ is also relative to $Z'$.

Theorem \ref{3th1} has a natural formulation for the case of
Diff-algebras relative to $Z'$. Given an ideal $I\subset \calo_Z$, and a smooth morphism $Z\to
Z'$, we define an extension of ideals $I \subset Diff^r_{Z'}(I)$,
$$Diff^r_{Z'}(I)(U)=\langle D(f) / f\in I(U) , D \in
Diff^r_{Z'}(U) \rangle $$ for each open $U$ in $Z$.

Since $Diff^r_{Z'}\subset Diff^{r+1}_{Z'}$, clearly
$Diff^r_{Z'}(I)\subset Diff^{r+1}_{Z'}(I)$ for $r\geq 0$.

Note finally that a Rees algebra $\bigoplus I_kW^k$
 over $Z$ (\ref{par3.0}) is a Diff-algebra relative to $Z'$, if
 and only if, for any positive integers $r\leq n$,
 $Diff^r_{Z'}(I_n)\subset I_{n-r}.$ In particular, for $Z'=Spec(k)$, condition ii) in Def \ref{3def1}
can be reformulated as:

ii') $Diff^r_{k}(I_n)\subset I_{n-r}.$ }

\end{parrafo}

\section{ Diff-algebras and singular locus.}

\begin{parrafo}

{\rm The notion Diff-algebra relative to a perfect field $k$, on a
smooth $k$-scheme $Z$, is closely related to the notion of {\em
order} at the local regular rings of $Z$. Recall that the order of
a non-zero ideal $I$ at a local regular ring $(R,M)$ is the
biggest integer $b$ such that $I\subset M^b$.

If $I\subset \calo_Z$ is a sheaf of ideals, $V(Diff^{b-1}_{k}(I))$
is the closed set of points of $Z$ where the ideal has order at
least $b$. We analyze this fact locally at a closed point $x$.

Let $\{x_1,\dots ,x_n\}$ be a regular system of parameters for
$\calo_{Z,x}$, and consider the differential operators
$\Delta^{\alpha}$, defined on $\calo_{Z,x}$ in terms of these
parameters, as in the Theorem \ref{3th1}. So at $x$,
$$(Diff^{b-1}_{k}(I))_x=\langle \Delta^{\alpha}(f)/ f\in I, 0\leq
| \alpha| \leq b-1\rangle.$$ One can now check at $\calo_{Z,x}$,
or at the ring of formal power series $\hat{\calo}_{Z,x}$, that
$Diff^{b-1}_{k}(I)$ is a proper ideal if and only if $I$ has order
at least $b$ at the local ring.

The operators $\Delta^{\alpha}$ are defined globally at a suitable
neighborhood $U$ of $x$. So if $\bigoplus I_n W^n \subset
O_{Z}[W]$ is a Diff-algebra relative to the field $k$ and $x\in
Z$ is a closed point, the Diff-algebra defined by localization, say $\bigoplus (I_n)_x\
W^n \subset O_{Z,x}[W]$, is properly included in
$\mathcal{O}_{Z,x}[W]$, if and only, for each index $k\in
\mathbb{N}$, the ideal $(I_k)_x$ has order at least $k$ at the
local regular ring $\mathcal{O}_{Z,x}$.}

\end{parrafo}

\begin{definition}\label{3defsing} The {\em singular locus} of a Rees algebra 
$\mathcal{G}=\bigoplus I_n W^n
\subset O_{Z}[W]$, will be  $$Sing(\mathcal{G})=\cap_{r\geq 0}
V(Diff^{r-1}_{k}(I_r)) (\subset Z).$$ It is the set of points
$x\in Z$ for which all $(I_r)_x$ have order at least $r$ (at
$\mathcal{O}_{Z,x}$).
\end{definition}

\begin{remark}\label{3rkdelsing}

 Assume that $f\in (I_r)_x$ has order $r$ at $\mathcal{O}_{Z,x}$, where $x$ is in $Sing(\mathcal{G})$.
Then, locally at $x$,  $Sing(\mathcal{G})$ is included in the set
of points of multiplicity $r$ (or say, $r$-fold points) of the
hypersurface $V(\langle f \rangle)$.

In fact $Diff^{r-1}_{k}(f)\subset Diff^{r-1}_{k}(I_r)$, and the
closed set defined by the first ideal is that of points of
multiplicity $r$.
\end{remark}

\begin{proposition}\label{3propsing}

1) If $\mathcal{G}= \bigoplus I_nW^n$ and $ \mathcal{G}'= \bigoplus I'_nW^n$ are Rees algebras with the
same integral closure (e.g. if $\mathcal{G}\subset \mathcal{G}'$
is a finite extension), then
$$Sing(\mathcal{G})=Sing(\mathcal{G}').$$

2) If $\mathcal{G}$  is a Rees algebra generated over $\calo_Z$ by
$\mathcal{F}=\{ g_{n_i}W^{n_i}, n_i>0 , 1\leq i\leq m \}$, then
$$Sing(\mathcal{G})=\cap V( Diff^{n_i}(\langle g_i \rangle ) ).$$

3) Let $\mathcal{G}''=\bigoplus I''_n W^n$ be the extension
of $\mathcal{G}$ to a Diff-algebra relative to $k$, as defined
in Theorem \ref{3th1}, then
$Sing(\mathcal{G})=Sing(\mathcal{G}'')$.

4) For any Diff-algebra $\mathcal{G}''=\bigoplus I''_n
W^n$, $Sing(\mathcal{G}'')=V(I''_1)$.

5) Let  $\mathcal{G}''=\bigoplus I''_n\cdot W^n$ be a
Diff-algebra. For any positive integer $r$,
$Sing(\mathcal{G}'')=V(I''_r)$.
\end{proposition}
\proof

1) The argument in \ref{IIpar3} shows that there is an index $N$,
so that $\mathcal{G}$ is finite over the subring $\bigoplus
I^k_{N}W^{Nk}$, and $\mathcal{G}'$ is finite over  $\bigoplus
I'^k_{N}W^{Nk}$. And furthermore, $I_N$ and $I'_N$ have the same
integral closure. In these conditions $Sing(\mathcal{G})$ is the
set of points $x\in Z$ where $I_N$ has order at least $N$ at
$\mathcal{O}_{Z,x}$, and similarly, $Sing(\mathcal{G}')$ is the
set of points $x\in Z$ where $I'_N$ has order at least $N$.
Finally, the claim follows from the fact that the order of an
ideal, at a local regular ring, is the same as the order of its
integral closure (\cite{ZS}, Appendix 3).

2) We have formulated 2) with a global condition on $Z$, however
this is always fulfilled locally. In fact, there is a covering of
$Z$ by affine open sets, so that the restriction of $\mathcal{G}$
is generated by finitely many elements. Let $U$ be such open set,
so $\mathcal{G}(U)=\bigoplus I_k(U)\cdot W^k$ is generated by
$\mathcal{F}=\{ g_{n_i}W^{n_i}, n_i>0 , 1\leq i\leq m \}$,
$g_{n_i}\in \mathcal{O}(U)$.

The claim is that
 $y\in Sing(\mathcal{G})\cap U$ if and only if
the order of $g_{n_i}$ at $\calo_{Z,y}$ is at least $n_i$, for
$1\leq i\leq m$.

The condition is clearly necessary. Conversely, if
$\mathcal{G}=\bigoplus I_n =\calo_U[\{g_iW^{n_i}\}_{g_iW^{n_i}\in
\mathcal{F}}]$, and each $g_{n_i}$ has order at at least $n_i$ at
$\calo_{Z,y}$, then $I_n$ (generated by weighted homogeneous
expressions on the $g_i$'s) has order at least $n$ at
$\calo_{Z,y}$.

3) We argue as in 2). Fix a regular system of parameters $\{x_1,\dots ,x_n\}$ at $x\in U$,
and differential operators $\Delta^{\alpha}$ as in the Theorem
\ref{3th1}. After suitable restriction we may assume that these operators are defined globally at $U$.

Formula (\ref{eq3422}) shows tht the Diff-algebra $\mathcal{G}''$, in the Theorem \ref{3th1}, is a
finite extension of the Rees algebra defined by
$$\mathcal{F'}=\{\Delta^{\alpha}(g_n)W^{n_i-\alpha}/
g_{n_i}W^{n_i}\in \mathcal{F}, \alpha=(\alpha_1, \alpha_2,\dots ,
\alpha_n) \in (\mathbb{N})^n, \mbox{ and } 0\leq |\alpha| <
n_i\}.$$

Note finally that if the order of $g_{n_i}$ at a local ring is
$\geq n_i$, then the order of $\Delta^{\alpha}(g_n)$ is $\geq
n_i-|\alpha|$.

4) The inclusion $Sing(\mathcal{G}'')\subset V(I''_1)$ holds, by
definition, for any Rees algebra. On the other hand, the
hypothesis ensures that $Diff^{r-1}(I''_r)\subset I''_1$, so
$Sing(\mathcal{G}'')\supset V(I''_1)$.

5) Follows from 4).
\endproof


\section{ On restrictions of Diff-algebras.}

The concept of Diff-algebra is defined here for Rees algebras over a scheme, say $V$, which is smooth over a field $k$. Let $V'$ be another smooth scheme over $k$ and let $V'\to V$ be a morphism of $k$-schemes, then there is a natural lifting of a Rees algebra over $V$ to a Rees algebra over $V'$. The goal in this section is Theorem \ref{th54} which states that the lifting of a Diff-algebra is again a Diff-algebra.

\begin{proposition}\label{3prop1}

 Let $V$ be smooth over a perfect field, and let $\mathcal{G}=\bigoplus I_k\cdot W^k$ be a Diff-algebra defined by ideals $I_k\subset \calo_{V}$.

\bigskip

A) If $V'\subset V$ is a closed and smooth subscheme, the restriction of
$\mathcal{G}$ to $V'$, say $$\mathcal{G}'=\bigoplus I_k
\calo_{V'}\cdot W^k,$$ is a Diff-algebra over $V'$.

\bigskip

B) If $V''\to V$ is a smooth morphism, then the natural extension,
say $$\mathcal{G}''=\bigoplus I_k\calo_{V''}\cdot W^k,$$ is a
Diff-algebra over $V''$.

\end{proposition}

\proof It is clear that both $\mathcal{G}'$ and $\mathcal{G}''$
are Rees algebras (\ref{par3.0}). We will show that conditions (i)
and (ii) in Definition \ref{3def1} hold.

It suffices to prove both results locally at closed points, say
$x\in Sing(\mathcal{G})$. Set $\mathcal{G}_x=\bigoplus I_k\cdot
W^k$ where now each $I_n$ is an ideal in $\calo_{V,x}$. We may
also replace the local ring by its completion.

A) Fix a closed point $x\in V'\subset V$ and a local regular
system of parameters, say $$\{x_1,\dots ,x_h,x_{h+1},\dots x_d\}$$
at $\calo_{V,x}$, such that $V'$ is locally defined by the ideal
$<x_1,\dots, x_h>$. Set $$\hat{\calo}_{V,x}=k'[[x_1,\dots
,x_h,x_{h+1},\dots x_d]],$$ where $k'$ is a finite extension of
$k$. For each multi-index $\alpha=(\alpha_1,\dots , \alpha_d)\in
\mathbb{N}^d,$
$$\Delta^{\alpha}=\Delta^{\alpha^{(1)}}\Delta^{\alpha^{(2)}};$$
where $\alpha^{(1)}=(\alpha_1,\dots , \alpha_h)\in \mathbb{N}^h$,
and $\alpha^{(2)}=(\alpha_{h+1},\dots , \alpha_d)\in
\mathbb{N}^{d-h}$.

Express an element $f_n\in I_n$ as
$$ f_n=\sum_{\alpha^{(1)} \in (\mathbb{N})^h}x_1^{\alpha_1}\cdots
x_h^{\alpha_h}a^{}_{\alpha^{(1)}},$$ $a^{}_{\alpha^{(1)}}\in
k'[[x_{h+1},\dots x_d]]$.

If $|\alpha^{(1)}|=\alpha_1+\dots +\alpha_h \leq n$, then
$a_{\alpha^{(1)}}W^{n-|\alpha^{(1)}|}$ is the class of
$\Delta^{\alpha^(1)}(f_n)W^{n-|\alpha^{(1)}|}$ in
$\hat{\calo}_{V',x}[W]$. So it is an element in the restricted
algebra. Similarly, if $|\alpha^{(1)}|+|\alpha^{(2)}| \leq n$,
$$\Delta^{\alpha^{(2)}}a_{\alpha^{(1)}}W^{n-|\alpha^{(1)}|-
|\alpha^{(2)}|}$$ is the class of the element
$(\Delta^{\alpha^{(2)}}\Delta^{\alpha^(1)})(f_n)\cdot
W^{n-|\alpha^{(1)}|-|\alpha^{(2)}|}$ in $\hat{\calo}_{V',x}[W].$

 For each index $m\geq 1$, $I_m\calo_{V'}\cdot W^m
$ is defined by the coefficient $a_{0}W^m$ ($0\in
(\mathbb{N})^h$), for each $f_mW^m \in I_m W^m$. Conditions (i)
and (ii) in \ref{3def1} are now easy to check.

For our further discussion we point out that $I_m\calo_{V'}
W^m $ also contains all coefficients
$a_{\alpha^{(1)}}W^{n-|\alpha^{(1)}|}$ of $f W^n \in I_nW^n$, and
$n-|\alpha^{(1)}|=m$.

B) Fix a point $x'\in V''$ mapping to $x\in V$. The completion of
$\calo_{V'',x'}$ contains that of $\calo_{V,x}$, say
$$\hat{\calo}_{V,x}=k'[[x_1,\dots ,x_d]]\subset \hat{\calo}_{V'',x'}=k'[[x_1,\dots
,x_d,x_{d+1},\dots x_e]].$$ Each ideal $I_n$ in $k'[[x_1,\dots
,x_d]]$ extends to $I_n\cdot k'[[x_1,\dots ,x_d,x_{d+1},\dots
x_e]]$; and the claim is that the extended algebra is a
Diff-algebra. The statement follows easily in this case, for
example by formula (\ref{eq3422}), which expresses generators of
the Diff-algebra in terms of generators of the Rees algebra.


\begin{definition}\label{3def07}

 Fix $\mathcal{G}=\bigoplus I_k\cdot W^k$, a Rees algebra over
$V$, and a morphism of $k$-schemes, say $V\stackrel{\pi}{\longleftarrow}  V'$. Assume that $V$ and $V'$ are smooth.
Define the {\em total transform} of $\mathcal{G}$ to be
$$\pi^{-1}(\mathcal{G})=\bigoplus I_k\calo_{V'}\cdot W^k.$$ Namely the Rees algebra
defined by the total transforms of the ideals $I_n$, $n\geq 0$.
\end{definition}

Note that the {\em restriction} in A) and the {\em natural
extension} in B), are particular examples of total transforms.

Assume that $V$ is affine and that $\mathcal{F}=\{g_{N_1}W^{N_1},\dots ,g_{N_s}W^{N_s}\}$, generate  $\mathcal{G}$ . Then each 
$g_{N_i}$ defines a global section, say $\pi^*(g_{N_i})$ on $V'$, and we set, say $$\pi^*(\mathcal{F})=\{\pi^*(g_{N_1})W^{N_1},\dots ,\pi^*(g_{N_s})W^{N_s}\}$$ as elements in $\calo_{V'}[W] $. 

\begin{lemma}\label{lemolvidado} Let $\mathcal{G}=\bigoplus I_k\cdot W^k (\subset \calo_V[W])$ be
a {\em Rees algebra} generated by a finite set
$\mathcal{F}=\{g_{N_1}W^{N_1},\dots ,g_{N_s}W^{N_s}\}$, and let
$V\stackrel{\pi}{\longleftarrow}  V'$ be a morphism of smooth schemes. Then
$\pi^{-1}(\mathcal{G})$ is generated by $\pi^*(\mathcal{F})$.
\end{lemma}
\proof

 Since any element of $I_M$ is a weighted homogeneous polynomial
expressions of degree $M$, in elements of $\mathcal{F}$, the total transform of
the ideal is
 also generated by elements that are weighted homogeneous on $\pi^*(\mathcal{F})$.
\endproof
In particular:

A) the restriction of $\mathcal{G}$ to $ V'(\subset V)$ is
generated by $\pi^*(\mathcal{F})=\{\overline{g}_{N_1}W^{N_1},\dots
,\overline{g}_{N_s}W^{N_1}\}$, where each $\overline{g}_{N_i}$ is
the restriction of $g_{N_i}$ to $V'$.

B) If $V\stackrel{\pi}{\longleftarrow}  V''$ is a smooth morphism, the total transform of
$\bigoplus I_k\cdot W^k$ to $V''$ is generated by
$\pi^*(\mathcal{F})=\{g_{N_1}W^{N_1},\dots ,g_{N_s}W^{N_s}\}$.

\begin{theorem}\label{th54} Let $V'\stackrel{\pi}{\longrightarrow} V$ be a morphism of smooth
schemes, then:

i) if $\mathcal{G}$ is a Diff-algebra over $V$, the total
transform $\pi^{-1}(\mathcal{G})$ is a Diff-algebra over $V'$.

ii)  $Sing(\pi^{-1}(\mathcal{G}))=\pi^{-1}(Sing(\mathcal{G}))$.
\end{theorem}
\proof Since $V'\stackrel{\pi}{\longrightarrow} V$ is of finite
type, it can be expressed locally in the form $V'\subset
V''\stackrel{\beta}{\longrightarrow} V,$ where $\beta$ is smooth.
So Prop \ref{3prop1} proves (i).

Fix a closed point $x\in Sing(\pi^{-1}(\mathcal{G}))$. Since
$Sing(\mathcal{G})=V(I_n)$ for all $n\geq 1$ (\ref{3propsing}), it
follows that $\pi(x)\in Sing(\mathcal{G})$. On the other hand, if
$\pi(x)\in Sing(\mathcal{G})$, the order of $I_n$ is at least $n$
at $\calo_{V,\pi(x)}$, for each $n\geq 1$; so the same holds at
$\calo_{V',x}$. This proves (ii).

\section{ On Diff-algebras and integral closures.}

\begin{parrafo}\label{rk55n-1}
{\em  The aim in this section is, essentially, the proof of Main Theorem \ref{TH515}.

The proof will require a better understanding of the notions of restriction already studied in in the last section. In this previous discussion restrictions
where studied for a closed immersion of smooth schemes, say $Z\subset V$. Here we will consider, at least for the first results, a closed immersion together with a retraction, say $V\to Z$.

Given $Z\subset V$ as above, a local retraction at a point $x\in Z$ can always be defined in an \'etale neighborhood.

Here, given a Rees algebra $\mathcal{G}=\bigoplus I_k\cdot W^k
(\subset \calo_V[W])$ (over $V$), the retraction $V\to Z$ will allow us to define a new Rees algebra over $Z$, called the Coefficient 
algebra.
%
%



Fix $x\in \Sing (\mathcal{G}) \cap Z$. The retraction defines an inclusion $\calo_{Z,x}\subset \calo_{V,x}$. 
Extend a regular system of parameters of $\calo_{Z,x}$, say $\{x_1,\dots ,x_h\}$, to  a regular system of parameters,
say $\{x_1,\dots ,x_h,x_{h+1},\dots x_d\}$, of $\calo_{V,x}$. We may assume here that $I(Z)$ is $<x_{h+1}, \dots , x_d>$ at $\calo_{V,x}$.
The construction of the coefficient algebra will be addressed firsts at the completion of the local rings.
So $\hat{\calo}_{V,x}$ is a ring of formal power series, say
$k'[[x_1,\dots ,x_h,x_{h+1},\dots x_d]]$, and $\hat{\calo}_{V',x}$
is $k'[[x_{h+1},\dots x_d]]$. The (induced)  local retraction is defined by
the inclusion $k'[[x_{h+1},\dots x_d]]\subset k'[[x_1,\dots ,x_h,x_{h+1},\dots
x_d]]$.

Set, as usual, $\mathcal{G}_x=\bigoplus I_k\cdot W^k (\subset
\calo_{V,x}[W])$, which also extends to a Rees algebra over
 $\hat{\calo}_{V,x}$. Express an element $f_n\in I_n$ as
$$ f_n=\sum_{\alpha^{(1)} \in (\mathbb{N})^h}x_1^{\alpha_1}\cdots
x_h^{\alpha_h}a^{}_{\alpha^{(1)}}, \hskip 0.5cm
a^{}_{\alpha^{(1)}}\in k'[[x_{h+1},\dots x_d]].$$

For any such $f_nW^n$, consider the set
$\{a^{}_{\alpha^{(1)}}\cdot W^{n-|\alpha^{(1)}|}, 0\leq
|\alpha^{(1)}|<n \}$, which we call the {\em coefficients of}
$f_nW^n$. So the coefficients of $f_nW^n$ is a finite set, defined
in terms of a regular system of parameters, and the weight of each
coefficient depends on the index $n$.

\bigskip

{\em Claim:} As $f_n W^n$ varies on the Rees algebra
$\mathcal{G}_x$, the coefficients of $f_nW^n$ generate a Rees
algebra, say $\mathcal{\mbox{Coeff}(G)}_x$, in $k'[[x_{h+1},\dots
x_d]][W]$.

The claim here is that the graded algebra
$\mathcal{\mbox{Coeff}(G)}_x$ is a finitely generated subalgebra
of $k'[[x_{h+1},\dots x_d]][W]$.

\bigskip

Assume that $\mathcal{F}=\{g_{N_1}W^{N_1},\dots ,g_{N_s}W^{N_s}\}$
generate $\mathcal{G}_x$. Express, for $1\leq i \leq s$:

\begin{equation}\label{eq5511}
g_{N_i}=\sum_{\alpha \in (\mathbb{N})^h}x_1^{\alpha_1}\cdots
x_h^{\alpha_h}a^{(i)}_{\alpha} \hskip 0.5cm a_{\alpha}\in
k'[[x_{h+1},\dots x_d]].
\end{equation}

We search for a finite set of coefficients, that span
$\mathcal{\mbox{Coeff}(G)}_x$. A first candidate would be
\begin{equation}\label{eqmlfp}
\mathcal{F}'_1=\{a^{(i)}_{\alpha}W^{N_i-|\alpha| }/ 0\leq |\alpha|
< N_i, 1\leq i \leq s\}.
\end{equation}
The difficulty appears already if we consider the product of two elements in $\mathcal{F}$, say
$g_{N_i}W^{N_i}\cdot g_{N_j}W^{N_j}=f_n W^n$ ($n=N_i+N_j$); and a
coefficient, say $a^{}_{\alpha^{(1)}}·W^{n-|\alpha^{(1)}|}$, of
$f_n W^n$.

It follows from \ref{eq5511} that
$$a^{}_{\alpha^{(1)}}=\sum_{\beta+\delta=\alpha^{(1)}}a^{(i)}_{\beta}a^{(j)}_{\delta},
$$
for $\beta$, $\delta$, and $\alpha^{(1)}$ in $(\mathbb{N})^h$.
Note that the previous expression cannot be formulated in  the form
$$a^{}_{\alpha^{(1)}}W^{n-|\alpha^{(1)}
 |}=\sum_{\beta+\delta=\alpha^{(1)}}a^{(i)}_{\beta}W^{N_i-|\beta|}a^{(j)}_{\delta}W^{N_j-|\delta
 |} .$$
In fact, it can happen that $|\delta|\geq N_j$, and we only consider $W$ with positive exponents. In particular, the previous
expression of $a^{}_{\alpha^{(1)}}W^{n-|\alpha^{(1)}|}$ is not
weighted homogeneous in $\mathcal{F}'_1$, and hence not in the
graded sub-algebra of $k'[[x_{h+1},\dots x_d]][W]$ generated by
$\mathcal{F}'_1$.

One way to remedy this situation is to allow $a^{(i)}_{\beta}$ to
have weight $n-|\alpha^{(1)}|$ if $|\delta
 |\geq N_j$. Note that in such case $$n-|\alpha^{(1)}|=N_i-|\beta|+N_j-|\delta|\leq N_i-|\beta|.$$
Therefore  $\mathcal{F}'_1$ can be enlarged to say,
\begin{equation}\label{eq5512}
\mathcal{F}_1=\{a^{(i)}_{\alpha}W^{n_{i,\alpha} }/ 0\leq |\alpha|
< N_i, 1\leq i \leq s, 0< n_{i,\alpha} \leq N_i-|\alpha|\},
\end{equation}
for $N_i$ and $\alpha$ as in $\mathcal{F}_1$; and the coefficients of $f_nW^n$ are now weighted homogeneous on
$\mathcal{F}_1$ (i.e., are in the sub-algebra of $k'[[x_{h+1},\dots
x_d]][W]$ generated by $\mathcal{F}_1$).

The argument applied here to $g_{N_i}W^{N_i}\cdot g_{N_j}W^{N_j}$,
also holds for the coefficients of any product of elements in
$\mathcal{F}$, and hence for the coefficients of any homogeneous
element in the algebra generated by
$\mathcal{F}=\{g_{N_1}W^{N_1},\dots ,g_{N_s}W^{N_s}\}$ (i.e., for
the coefficients of any homogeneous element of $\mathcal{G}_x$).

This shows that there is an inclusion of subalgebras in
$k'[[x_{h+1},\dots x_d]][W]$, say
\begin{equation}\label{eq5535}k'[[x_{h+1},\dots x_d]][\mathcal{F}'_1] \subset
\mathcal{\mbox{Coeff}(G)}_x\subset k'[[x_{h+1},\dots
x_d]][\mathcal{F}_1].
\end{equation}
On the other hand $k'[[x_{h+1},\dots x_d]][\mathcal{F}'_1]\subset
k'[[x_{h+1},\dots x_d]][\mathcal{F}_1]$ is a finite extension
(\ref{IIrk1}, 2)). In particular $\mathcal{\mbox{Coeff}(G)}_x$ is
finitely generated. }
\end{parrafo}

\begin{remark}\label{rk55m}
1) $\mathcal{F}'_1$ can be extended to a finite set, say
$\mathcal{F}''_1$, of generators of $\mathcal{\mbox{Coeff}(G)}_x$.

2) $Sing(k'[[x_{h+1},\dots x_d]][\mathcal{F}'_1])=
Sing(\mathcal{\mbox{Coeff}(G)}_x)$ (Prop \ref{3propsing}, (1)).

3) $Sing(\mathcal{\mbox{Coeff}(G)}_x)$ can be naturally identified
with the intersection $Z\cap Sing(\mathcal{G})$ locally at the
point $x$.

To check this last point 3) note first that the singular locus of
$k'[[x_{h+1},\dots x_d]][\mathcal{F}'_1]$ can be naturally
identified with the intersection $Z\cap Sing(\mathcal{G})$. This
follows from the definition of $\mathcal{F}'_1$ in (\ref{eqmlfp}),
and the expressions in (\ref{eq5511}). Finally apply 2).

\end{remark}

\begin{parrafo}\label{rk55n} {\em Fix, as in \ref{rk55n-1}, an inclusion of smooth schemes
$Z \subset V$, and a retraction say $V\to Z$, defined locally at a closed point $x \in Z$.  Let $x_1,\dots, x_h$ be
part of a regular system of parameters for $\calo_{V,x}$ so that
$<x_1,\dots, x_h>$ defines $I(Z)$ at $\calo_{V,x}$; and let
$\{x_{h+1},\dots x_d\}$ be a regular system of parameters for
$\calo_{Z,x}$. 
%
%

Given $\mathcal{G}=\bigoplus I_k\cdot W^k (\subset \calo_V[W])$,
we have defined $\mathcal{\mbox{Coeff}(G)}$ at
${\hat{\calo}}_{Z,x}[W]$. We now show that it can also be defined
in ${{\calo}}_{Z,x}[W]$, and that the definition relies on the inclusion $Z\subset V$ and on the retraction. Express an element $f_n\in I_n{\hat{\calo}}_{Z,x}$ as
$$ f_n=\sum_{\alpha^{(1)} \in (\mathbb{N})^h}x_1^{\alpha_1}\cdots
x_h^{\alpha_h}a^{}_{\alpha^{(1)}},$$ $a^{}_{\alpha^{(1)}}\in
k'[[x_{h+1},\dots x_d]]$. For each multi-index $\alpha^{(1)}$, $0
\leq |\alpha^{(1)}|\leq n$, the coefficient $a^{}_{\alpha^{(1)}}$
can be identified with the class of $\Delta^{\alpha^{(1)}}(f_n)$
in ${\hat{\calo}}_{Z,x}$. However,  $\Delta^{\alpha^{(1)}}$ is a
differential operator, relative to the local retraction $V \to Z$,
 $\Delta^{\alpha^{(1)}}(f_n)$
in an element in ${{\calo}}_{V,x}$, and one can consider
the class of this element in ${{\calo}}_{Z,x}$.

This shows that $\mathcal{\mbox{Coeff}(G)}(\subset
{{\calo}}_{Z,x}[W])$, is the restriction via $Z\subset V$, of the
extension of $\mathcal{G}$ defined by the Diff-algebra relative
to the local retraction (\ref{par3rela}). In other words, $\mathcal{\mbox{Coeff}(G)}(\subset
{{\calo}}_{Z,x}[W])$ is defined in terms of:

i) the surjection $\calo_{V,x} \to \calo_{Z,x}$; and

ii) the inclusion $ \calo_{Z,x}\subset \calo_{V,x}$ (\'etale
locally).

This proves the following Remark in the formal case:}

\end{parrafo}
\begin{remark}\label{64}
Consider $\mathcal{G}\subset {\calo}_{V,x}[W]$, and, as before, a
surjection $\hat{\calo}_{V,x} \to \hat{\calo}_{Z,x}$; and an
inclusion $\hat{\calo}_{Z,x}\subset \hat{\calo}_{V,x}$. Let
$\{x_{h+1},\dots ,x_d\}$ be a regular system of parameters for
${\calo}_{Z,x}$, and assume that $\{x_1, \dots x_h,x_{h+1},\dots
,x_d\}$ and $\{x'_1, \dots x'_h,x_{h+1},\dots , x_d\}$ are two
extensions to regular system of parameters for ${\calo}_{V,x}$, and
that $<x_1, \dots ,x_h>=<x'_1, \dots ,x'_h>=I(Z)\subset
{\calo}_{V,x}$

The same inclusion $ \calo_{Z,x}\subset \calo_{V,x}$ can be
expressed as
$$\hat{\calo}_{Z,x}=k'[[x_{h+1},\dots , x_d]]\subset k'[[x_1,\dots
,x_h,x_{h+1},\dots , x_d]]=\hat{\calo}_{V,x}$$ or as
$$\hat{\calo}_{Z,x}=k'[[x_{h+1},\dots , x_d]]\subset k'[[x'_1,\dots
,x'_h,x_{h+1},\dots , x_d]]=\hat{\calo}_{V,x}.$$

In both cases  $Coeff(\mathcal{G})\subset \hat{\calo}_{Z,x}[W]$ is
the same.

The discussion in \ref{rk55n} also shows that, of course, the
definition of $Coeff(\mathcal{G})\subset \hat{\calo}_{Z,x}[W]$ is
independent of the coordinates we choose in the subring
$\hat{\calo}_{Z,x}$.
\end{remark}

\begin{parrafo} {\em Set $\mathcal{G}\subset {\calo}_{V,x}[W]$ and
$\hat{\calo}_{Z,x}=k'[[x_{h+1},\dots , x_d]]\subset k'[[x_1,\dots
,x_h,x_{h+1},\dots , x_d]]=\hat{\calo}_{V,x}$ as above, where
$<x_1, \dots ,x_h>=I(Z)\subset \hat{\calo}_{V,x}$. Note that the
natural identification of $k'[[x_{h+1},\dots , x_d]]$ with
$\hat{\calo}_{V,x}/I(Z)$, provides an inclusion:
$$\overline{\mathcal{G}}\subset Coeff(\mathcal{G})\subset
\hat{\calo}_{Z,x}[W],$$ where $\overline{\mathcal{G}}\subset
(\hat{\calo}_{V,x}/I(Z))[W]$ denotes the restriction. Furthermore,
this inclusion is an equality if $\mathcal{G}$ is a Diff-algebra:

}
\end{parrafo}

\begin{lemma} \label{lem5885} With the setting as above, the restriction of
$G(\mathcal{G})$ to the smooth subscheme $Z$ is
 the Diff-algebra spanned by
$\mathcal{\mbox{Coeff}(G)}$ (i.e., the Diff-algebra generated by
$\mathcal{\mbox{Coeff}(G)}$ in ${\calo}_{Z,x}[W]$).
\end{lemma}
\proof
 The previous discussion shows that
$\mathcal{\mbox{Coeff}(G)}$ is included in the restriction of
$G(\mathcal{G})$, which is a Diff-algebra over ${\calo}_{Z,x}$
(Prop \ref{3prop1},A)). In particular, the Diff-algebra spanned
by $\mbox{Coeff}(G)_x$ is included in the restriction. The claim
is that this last inclusion is an equality.

Here $G(\mathcal{G})=\bigoplus I'_k\cdot W^k$ is the
Diff-algebra generated by $\mathcal{G}$, so to prove this
equality it suffices to show that given $f_n\in I_n$, and
$\alpha=(\alpha_1,\dots , \alpha_d)\in (\mathbb{N})^d,$ $0\leq
|\alpha|<n$, the class of $\Delta^{\alpha}(f_n)W^{n-|\alpha|}$ in
${{\calo}}_{Z,x}[W]$, is in the Diff-algebra generated by
$\mathcal{\mbox{Coeff}(G)}$.

For this last claim
 we argue as in the proof of Prop
\ref{3prop1}, (A), by splitting each multi-index
$\alpha=(\alpha_1,\dots , \alpha_d)\in (\mathbb{N})^d:$
$$\Delta^{\alpha}=\Delta^{\alpha^{(1)}}\Delta^{\alpha^{(2)}};$$
where $\alpha^{(1)}=(\alpha_1,\dots , \alpha_h)$, and
$\alpha^{(2)}=(\alpha_{h+1},\dots , \alpha_d)$.

The class of $\Delta^{\alpha^{(1)}}(f_n)W^{n-|\alpha^{(1)}|}$ is
$a^{}_{\alpha^{(1)}}W^{n-|\alpha^{(1)}|}\in
\mathcal{\mbox{Coeff}(G)}$; and that of
$\Delta^{\alpha}(f_n)W^{n-|\alpha|}$ is
$\Delta^{\alpha^{(2)}}(a_{\alpha^{(1)}})W^{n-|\alpha^{(1)}|-|\alpha^{(2)}|}$,
which is clearly in the Diff-algebra spanned by
$\mathcal{\mbox{Coeff}(G)}$.

\begin{corollary}\label{clm5885} Fix a smooth scheme $V$, a Rees algebra $\mathcal{G}$, and a closed and smooth subscheme
$Z$ of $V$. If $G(\mathcal{G})$ denotes the Diff-algebra spanned
by $\mathcal{G}$, and if $[G(\mathcal{G})]_Z$ denotes the
restriction to $Z$, then $Sing([G(\mathcal{G})]_Z)$, as closed set
in $Z$, can be identified with $Z\cap Sing(\mathcal{G})$.
\end{corollary}
This follows from Lemma \ref{lem5885} and \ref{rk55m}, 3). In fact, a local retraction can  be defined at an \'etale neighborhood of a closed point $x\in Z\cap Sing(\mathcal{G})$, and then argue at the local rings. 

\begin{remark}\label{rk56}{\em On the one dimensional case.}

We discuss here some particular features of the $G$ operator on Rees algebras, which hold when the dimension of the underlying smooth scheme is one. Let $\mathcal{G}=\bigoplus I_k\cdot W^k (\subset \calo_{V'}[W])$
be a {\em Rees algebra} over a one dimensional smooth scheme $V'$.
The aim is to prove that in the one-dimensional case $\mathcal{G}\subset G(\mathcal{G})$ is a finite extension.

If we assume that some $I_k\neq 0$, then $Sing(\mathcal{G})$ is a
finite set of points.
Fix $x\in Sing(\mathcal{G})$ and set $\hat{\calo}_{V',x}=k'[[t]]$,
so
$$\mathcal{G}=\bigoplus_{r\geq 1} <t^{a_r}>\cdot W^r,$$
and $a_r \geq r$ for each index $r$. Define
$$\lambda_{\mathcal{G}} =inf_r\{ \frac{a_r}{r}\},$$
and note that $\lambda_{\mathcal{G}} \geq 1$.

Let $\{g_{N_1}W^{N_1},\dots ,g_{N_s}W^{N_s}\}$ be a set of
generator locally at a closed point $x\in Sing(\mathcal{G})$.
 Fix any integer $M$ divisible by all $N_i$, $1\leq i \leq s$,
then $$\lambda_{\mathcal{G}}= \frac{\nu(I_M)}{M}$$ where
$\nu(I_M)$ denotes the order of the ideal at $\calo_{V',x}$. Let
$\widetilde{\mathcal{G}}$ denote the integral closure of
${\mathcal{G}}$.

{\em Claim 1:} The integral closure of ${\mathcal{G}}$ is
determined by the rational number $\lambda_{\mathcal{G}}$, and
$\lambda_{\mathcal{G}}=\lambda_{\widetilde{\mathcal{G}}}$.

In fact, by usual arguments of toric geometry, we conclude that
$t^n\cdot W^m \in \widetilde{\mathcal{G}}$, if and only if
$\frac{n}{m} \geq \lambda_{\mathcal{G}}$. This proves the claim.

\bigskip

Recall that $Sing(\mathcal{G})=Sing(G(\mathcal{G}))$.

{\em Claim 2:} Locally at any $x\in Sing(\mathcal{G})$, both
$\mathcal{G}$ and $G(\mathcal{G})$ have the same integral closure.

Let $\Delta^r$, $r\geq 0$, be defined in terms of the Taylor development in $k'[[t]]$, as in the proof of Theorem \ref{3th1}.
We prove our claim by showing that
$\lambda_{\mathcal{G}}=\lambda_{G(\mathcal{G})}$. To this end note
that given $t^a\cdot W^b \in \mathcal{G}$, and an operator
$\Delta^r$, $0 \leq r < b$,  $$\Delta^r(t^a)\cdot W^{b-r}=d\cdot
t^{a-r}\cdot W^{b-r},$$ where $d$ is the class of an integer in
the field $k'$. Since $a\geq b
>r\geq 0$ it follows that
$\frac{a-r}{b-r}\geq \frac{a}{b}$, so Claim 2 follows from Claim
1.
\end{remark}



\begin{parrafo}\label{par5144}{\em

Let $\mathcal{G}$ be a Rees algebra over $V$, and assume, after  restriction to affine open set, that it is
generated by $\{g_{N_1}W^{N_1},\dots ,g_{N_s}W^{N_s}\}$. Let $M$ is a positive
integer divisible by all $N_j$, $1\leq j \leq s$; and consider the
Rees ring $\calo_V [I_MW^M]$. Recall that $\calo_V [I_MW^M]\subset
\mathcal{G}$ is a finite extension of graded algebras, and
that any Rees algebra is a finite extension of a Rees ring of
an ideal (\ref{IIpar3}).

Given two Rees algebras $\mathcal{G}_1=\bigoplus_{r\geq 0}
I(1)_rW^{r}$ and $\mathcal{G}_2=\bigoplus_{r\geq 0} I(2)_rW^{r}$,
there is always a positive integer $M$ such that both are integral
extensions of the Rees ring generated by the $M$-th term, say
$\bigoplus_{k\geq 0} I(1)^k_MW^{km}$ and $\bigoplus_{k\geq 0}
I(2)^k_MW^{km}$.}

\end{parrafo}

\begin{proposition}\label{prop58}
Fix two Rees algebras $\mathcal{G}_1$ and $\mathcal{G}_2$ over a
smooth scheme $V$ over a field $k$. Assume that for any morphism
of regular $k$-schemes, say $V'\stackrel{\pi}{\longrightarrow} V$
, where $V'$ is one dimensional,
both pull-backs have the same integral closure (i.e., that
 $\widetilde{\pi^{-1}(\mathcal{G}_1)}=\widetilde{\pi^{-1}(\mathcal{G}_2)}$).
Then $\mathcal{G}_1$ and $\mathcal{G}_2$ have the same integral
closure in $V$.
\end{proposition}
\proof Fix a positive integer $M$ and ideals ideals $I(1)_M$ and $I(2)_M$ as in
\ref{par5144}. We may assume here that $\pi$ is of finite type.
Lemma \ref{lemolvidado} and the previous properties show that under the condition of the
hypothesis both $I(1)_M$ and $I(2)_M$ have the same integral
closure in $\calo_V$ (\ref{lemclauent}). In particular,
$\mathcal{G}_1$ and $\mathcal{G}_2$ have the same integral
closure.

\begin{proposition}\label{prop59}
Let $\mathcal{G}_1\subset \mathcal{G}_2 (\subset \calo_{V}[W])$ be
a finite extension of Rees algebras over a smooth scheme $V$, and
let $V'$ be a smooth one dimensional subscheme in $V$. Fix $x\in
V'$ and a regular system of coordinates $\{ x_1,\ldots
,x_{d-1},x_d\}$ for $\calo_{V,x}$, so that the curve is locally
defined by $<x_1,\ldots ,x_{d-1}>$. Then
$$\mbox{Coeff}(\mathcal{G}_1)\subset \mbox{Coeff}(\mathcal{G}_2)$$
is a finite extension in $\calo_{V'}[W]$.
\end{proposition}
\proof



By fixing coordinates we also fix a local retraction at an \'etale neighborhood of $x\in V$. Strictly speaking the coefficient algebras are defined in such neighborhood. We sometimes work at the completions
to ease the notation. Express any $ f\in \hat{\calo}_V=k'[[x_1,\ldots ,x_{d-1},x_d]]$
as:

\begin{equation*}
f=\sum_{\alpha \in (\mathbb{N})^{d-1}}x_1^{\alpha_1}\cdots
x_{d-1}^{\alpha_{d-1}}a^{}_{\alpha} \hskip 0.5cm a_{\alpha}\in
k'[[ x_d]].
\end{equation*}
The coefficients of $fW^N$ are $\{a_{\alpha}W^{N-|\alpha|} / 0
\leq |\alpha|< N\}$, and we define
$$sl_{V'}(fW^N)=min \{\frac{\nu(a^{}_{\alpha})}{N-|\alpha| }/ 0\leq |\alpha| < N\},$$
where $\nu(a^{}_{\alpha})$ denotes the order of $a^{}_{\alpha}$ in
$ k'[[ x_d]]$. Set $\mbox{Coeff}(\mathcal{G}_1)$ and $
\mbox{Coeff}(\mathcal{G}_2)$ in $\calo_{V'}[W]$, as in
(\ref{rk55n-1}). Assume that $\mathcal{F}_1=\{
f_{N_1}W^{N_1},\dots ,f_{N_s}W^{N_s}\}$ generate $\mathcal{G}_1$
locally at $x$, and that $\mathcal{F}_2=\{ g_{M_1}W^{M_1},\dots
,g_{M_t}W^{M_t}\}$ generate $\mathcal{G}_2$.
The inclusion $\mbox{Coeff}(\mathcal{G}_1)\subset
\mbox{Coeff}(\mathcal{G}_2)$ at $\calo_{V'}[W]$ is clear.

Set $\mbox{Coeff}(\mathcal{G}_i)=\bigoplus_{r\geq 0} J(i)_rW^{r}$
in $\calo_{V'}[W]$, for $i=1,2$. Note that $J(1)_r=0$ for all
$r\geq 1$ iff $V' \subset Sing(\mathcal{G}_1)$ iff
$\mathcal{G}_1\subset \bigoplus_{r\geq 0} P^{r}W^{r}$; where $P$
is the ideal defining the smooth subscheme $V'$. Since
$\mathcal{G}_1\subset \mathcal{G}_2$ is finite, it follows that
also $J(2)_r=0$ for all $r\geq 1$.

Assume now that some $J(1)_r$ is not zero for some $r>0$. The
inclusion $\mbox{Coeff}(\mathcal{G}_1)\subset
\mbox{Coeff}(\mathcal{G}_2)$ ensures that
\begin{equation}\label{eqlam1551}
\lambda_{\mbox{Coeff}(\mathcal{G}_1)}\geq
\lambda_{\mbox{Coeff}(\mathcal{G}_2)},
\end{equation}
and we shall prove the claim, in what follows, by showing that they are equal (see Remark \ref{rk56}).

Each  $g_{M_j}W^{M_j}$ is integral over the localization of
$\mathcal{G}_1$ in $\calo_{V,x}[W]$; and this property is
preserved by any change of rings. Namely, for any ring
homomorphism $\phi: \calo_{V,x}\to S$, $\phi(\mathcal{G}_1)$ is a
Rees algebra in $S[W]$, and $\phi(g_{M_j})W^{M_j}$ is integral
over $\phi(\mathcal{G}_1)$.

Express, for any $g_{M_j}W^{M_j}\in \mathcal{F}_2$:
\begin{equation}\label{eq55112}
g_{M_j}=\sum_{\alpha \in (\mathbb{N})^h}x_1^{\alpha_1}\cdots
x_h^{\alpha_h}a^{(j)}_{\alpha} \hskip 0.5cm a_{\alpha}\in k'[[
x_d]],
\end{equation}
and set
$$\mathcal{F}'_2=\{a^{(j)}_{\alpha}W^{M_j-|\alpha| }/ 0\leq |\alpha| < M_j,
1\leq j \leq t\}$$ (coefficients of all $g_{M_j}$'s).

We know that $k'[[x_d]][\mathcal{F}'_2] \subset
\mbox{Coeff}(\mathcal{G}_2)$ is a finite extension in
$k'[[x_d]][W]$ (see \ref{eq5535}); in particular:
$$\lambda_{\mbox{Coeff}(\mathcal{G}_2)}=min \{\frac{\nu(a^{(j)}_{\alpha})}{M_j-|\alpha| }/ 0\leq |\alpha| < M_j,
1\leq j \leq t\} \hskip 0.5cm (\ref{rk56}),$$ or, equivalently:
$$\lambda_{\mbox{Coeff}(\mathcal{G}_2)}=min \{sl_{V'}( g_{M_j}),
1\leq j \leq t\}.$$

So equality in (\ref{eqlam1551}) would follow if we show that
$\lambda_{\mbox{Coeff}(\mathcal{G}_1)} \leq
\frac{\nu(a^{(j)}_{\alpha})}{M_j-|\alpha| }$ for each fraction as
above. 

We will assume that
\begin{equation}\label{eqlapu1551}
\frac{\nu(a^{(j_0)}_{\alpha})}{M_{j_0}-|\alpha|} <
\lambda_{\mbox{Coeff}(\mathcal{G}_1)}
\end{equation} 
for some index $1\leq j_0
\leq t$, or equivalently, that $sl_{V'}(g_{M_{j_0}})<
\lambda_{\mbox{Coeff}(\mathcal{G}_1)}$ for some index $j_0$, and
show that in such case $g_{M_{j_0}}W^{M_{j_0}}$ is not integral
over $\mathcal{G}_1$; which is a contradiction.

Define, as before, $sl_{V'}(f_{N_i})$ for each $f_{N_i}W^{N_i}\in
\mathcal{F}_1$, so that
$$\lambda_{\mbox{Coeff}(\mathcal{G}_1)}=min\{sl_{V'}(f_{N_i});
1\leq i\leq s\}.$$

We will show that if $sl_{V'}(g_{M_{j_0}})<
\lambda_{\mbox{Coeff}(\mathcal{G}_1)}$, for some index $j_0$, a
ring $S$ and a morphism $\phi: \hat{\calo}_V=k'[[x_1,\ldots
,x_{d-1},x_d]] \to S$ can be defined so that
$\phi(g_{M_{j_0}})W^{M_{j_0}}$ is not integral over
$\phi(\mathcal{G}_1)$.

Given $f=\sum_{\alpha} \lambda_{\alpha}x_1^{\alpha_1}\dots
x_{d-1}^{\alpha_{d-1}}x_d^{\alpha_d}\in k'[[x_1,\ldots
,x_{d-1},x_d]]$, set
$$Supp(f)=\{\alpha \in \mathbb{N}^d / \lambda_{\alpha}\neq 0\}.$$

Let $a>0$ and $b>0$ be positive integers such that
$$\lambda=\lambda_{\mbox{Coeff}(\mathcal{G}_1)}=\frac{a}{b}.$$
Define $l:\mathbb{R}^d\to \mathbb{R}$, $l(y_1,\dots
,y_d)=ay_1+ay_2+\cdots ay_{d-1}+by_d$, which maps $\mathbb{N}^d$
into $\mathbb{N}$.

It follows that for a fixed integer $N$:
$$l(N,0, \dots, 0,0)=l(0,  N,\dots,0,0) =\cdots =l(0,\dots,N, 0)= l(0,\dots,0,\lambda N)= aN.$$
Given $(\alpha_1, \dots, \alpha_{d-1},s)\in \mathbb{N}^{d}$, if
$l(\alpha_1,\alpha_2, \dots, \alpha_{d-1},s)< aN$,
$|\alpha|:=\alpha_1+ \dots +\alpha_{d-1}<N$. Furthermore:
\begin{equation}\label{eq61}
l(\alpha_1,\alpha_2, \dots, \alpha_{d-1},s)< aN \Leftrightarrow
a|\alpha|+bs<aN \Leftrightarrow \frac{s}{N-|\alpha|}<\lambda.
\end{equation}

We show now that:

{\bf 1)} For each $f_{N_i}W^{N_i}\in \mathcal{F}_1$,
$Supp(f_{N_i})$ is included in the half space $l(y_1,\dots
,y_d)\geq aN_i.$

{\bf 2)} For some $f_{N_i}W^{N_i}\in \mathcal{F}_1$, the
intersection of $Supp(f_{N_i})$ with the hyperplane $l(y_1,\dots
,y_d)= aN_i$ is not empty.

{\bf 3)} For some $g_{M_{j_0}}W^{M_{j_0}}\in \mathcal{F}_2$,
$Supp(g_{M_{j_0}})$ is not included in the half space $l(y_1,\dots
,y_d)\geq aM_{j_0}.$

\bigskip

In order to prove 1) set
\begin{equation}\label{eqtico}
f_{N_i}=\sum_{\alpha \in (\mathbb{N})^{d-1}}x_1^{\alpha_1}\cdots
x_{d-1}^{\alpha_{d-1}}a^{(i)}_{\alpha} \hskip 0.5cm
a^{(i)}_{\alpha}\in k'[[ x_d]].
\end{equation}
and assume that
$$x_1^{\alpha_1}\cdots
x_{d-1}^{\alpha_{d-1}}x_d^s$$ is a monomial with non-zero
coefficient in this expression (i.e., assume that $(\alpha_1,\dots
,\alpha_{d-1},s)\in Supp(f_{N_i})$). The claim in 1) is that
$l(\alpha_1,\dots ,\alpha_{d-1},s)\geq aN_i$. In fact, if
$l(\alpha_1,\dots ,\alpha_{d-1},s)<aN_i$, then
$|\alpha|:=\alpha_1+ \dots +\alpha_{d-1}<N_i$ and
$\frac{s}{N_i-|\alpha|}<\lambda$ (\ref{eq61}). But in such case
$$sl_{V'}(f_{N_i})\leq \frac{s}{N_i-|\alpha|}< \lambda=\lambda_{\mbox{Coeff}(\mathcal{G}_1)}=min\{sl_{V'}(f_{N_i});
1\leq i\leq s\},$$ which is a contradiction.

Both conditions 2) and 3) follow similarly, from (\ref{eqlapu1551}) and ((\ref{eq61}).

\bigskip

 Set $S=k''[[t]]$ for some field extension $k''$
of $k'$, and define $\beta :k'[[x_1,\dots ,x_d]]\to k''[[t]]$ the
continuous morphism, such that $\beta(x_i)=\lambda_i t^a$
($\lambda_i\in k''$), for $1\leq i \leq d-1$, and
$\beta(x_d)=t^b$.

So $\beta(\mathcal{G}_1)$ is the Rees algebra in $k''[[t]][W]$
generated by $\{ \beta(f_{N_i})W^{N_i}, 1\leq i \leq s\}$.

It follows that for $k''$ an infinite field, and for sufficiently
general $\lambda_i\in k''$:

{\bf 1')} $\beta(f_{N_{i}})$ has order at least $aN_{i}$ in $k''[[t]]$.

{\bf 2')} $\beta(f_{N_{i_0}})$ has order $aN_{i_0}$ for some
$f_{N_{i_0}}W^{N_{i_0}}\in \mathcal{F}_1$.

{\bf 3')} $\beta(g_{M_{j_0}})$ has order strictly smaller then
$aM_{j_0}$.

Finally Claim 1 in Remark \ref{rk56}, where now $\lambda_{
\beta(\mathcal{G}_1)}=a$, asserts that
$\beta(g_{M_{j_0}})W^{M_{j_0}}$ is not integral over
$\beta(\mathcal{G}_1)$; so  (\ref{eqlapu1551}) can not hold.
\endproof

The following Theorem can also be proved by other means, which involve Hironaka's theory
on infinitely near points in \cite{Hironaka03}; a theory based on
the behavior by monoidal transforms. Our proof relies on the
previous development in this section, which will also be used for
the proof of Theorem \ref{TH513}.

\begin{theorem} (Main Theorem) \label{TH515} Let $\mathcal{G}_1\subset \mathcal{G}_2$ be an
inclusion of Rees algebras over a smooth scheme $V$. Let
$G(\mathcal{G}_i)$ be the Diff-algebra spanned by
$\mathcal{G}_i$ $(i=1,2)$ . If $\mathcal{G}_1\subset
\mathcal{G}_2$ is a finite extension, then
$G(\mathcal{G}_1)\subset G(\mathcal{G}_2)$ is a finite extension.
\end{theorem}
\proof

The inclusion $G(\mathcal{G}_1)\subset G(\mathcal{G}_2)$ is clear.
We will argue locally at a point $x\in Sing(\mathcal{G}_1)$, and
we make use of the criterion in Proposition \ref{prop58} to show
that the extension is finite. Let $\mathcal{F}_1=\{
f_{N_1}W^{N_1},\dots ,f_{N_s}W^{N_s}\}$ generate $\mathcal{G}_1$
locally at $x$, and let $\mathcal{F}_2=\{ g_{M_1}W^{M_1},\dots
,g_{M_t}W^{M_t}\}$ generate $\mathcal{G}_2$.

Set $\pi: V' \to V$ where $V'$ is one dimensional, and let $x'\in
V'$ map to $x$. Locally at $x'$ one can factor $\pi$ as
$V'\subset V''\to V$, so that $\phi: V''\to V$ is smooth.  Let
$\phi^{-1}(\mathcal{G}_1)$, $\phi^{-1}(\mathcal{G}_2)$ denote the
total transforms of $\mathcal{G}_1$, $\mathcal{G}_2$; and
$\phi^{-1}(G(\mathcal{G}_1))$, $\phi^{-1}(G(\mathcal{G}_2))$ be
the total transforms of $G(\mathcal{G}_1)$, $G(\mathcal{G}_2)$.

If $\{x_1,\dots ,x_d\}$ is a regular system of parameters for
$\calo_{V,x}$, then $\{x_1,\dots ,x_d\}$ extends to a regular
system of parameters, say $\{x_1,\dots ,x_d,\cdots ,x_e\}$ for
$\calo_{V'',x'}$. It is easy to check that

1) $\mathcal{F}_1=\{ f_{N_1}W^{N_1},\dots ,f_{N_s}W^{N_s}\}$
generate $\phi^{-1}(\mathcal{G}_1)$ locally at $\calo_{V'',x'}$;

2) $\mathcal{F}_2=\{ g_{M_1}W^{M_1},\dots ,g_{M_t}W^{M_t}\}$
generate $\phi^{-1}(\mathcal{G}_2)$ at $\calo_{V'',x'}$.

3) $\phi^{-1}(G(\mathcal{G}_1))$ is the Diff-algebra generated
by $\phi^{-1}(\mathcal{G}_1)$.

4) $\phi^{-1}(G(\mathcal{G}_2))$ is the Diff-algebra generated
by $\phi^{-1}(\mathcal{G}_2)$.

Therefore the setting at $V$ and at $V''$ is the same, and hence,
in order to apply Proposition \ref{prop58} we need only to show
that given a finite extension $\mathcal{G}_1\subset
\mathcal{G}_2$, the {\em restrictions} of the Diff-algebras
$G(\mathcal{G}_i)$, $i=1,2$, to a smooth one dimensional scheme
$V'$, have the same integral closure.

 Lemma \ref{lem5885} says that the restriction of $G(\mathcal{G}_i)$
to $V'$ is the Diff-algebra generated by
$\mbox{Coeff}(\mathcal{G}_i)$ ($i=1,2$). Remark \ref{rk56} shows
that for each index $i=1,2$, the Rees algebra
$\mbox{Coeff}(\mathcal{G}_i)$, and the Diff-algebra generated by
$\mbox{Coeff}(\mathcal{G}_i)$, have the same integral closure. So
it suffices to show that $\mbox{Coeff}(\mathcal{G}_1)$  and
$\mbox{Coeff}(\mathcal{G}_2)$ have the same integral closure,
which was proved in Prop \ref{prop59}. In fact, 1),2),3), and 4)
ensure that the setting of Prop \ref{prop59} hold.

\begin{theorem}\label{TH513} Let $\mathcal{G}_1\subset \mathcal{G}_2$ be an
inclusion of Rees algebras over a smooth scheme $V$. Fix a smooth
subscheme $Z\subset V$, and a local (or formal) retraction $ V\to
Z$. If $\mathcal{G}_1\subset \mathcal{G}_2$ is a finite extension,
then $\mbox{Coeff}(\mathcal{G}_1)\subset
\mbox{Coeff}(\mathcal{G}_2)$ is also finite.
\end{theorem}
\proof

Set $\pi: C \to Z$ where $C$ is smooth and one dimensional, and
let $x'\in C$ map to $x$. Locally at $x'$, one can factor $\pi$ as
$C\subset Z_1\to Z$, so that $\phi: Z_1\to Z$ is smooth. The
retraction of $V$ on $Z$, together with the morphism $Z_1\to Z$,
define by fiber products, a retraction say $V_1\to Z_1$, and a
smooth morphism, say $V_1\to V$.

The total transform of $\mathcal{G}_1\subset \mathcal{G}_2
(\subset \calo_V[W] )$ to say $\mathcal{G}'_1\subset
\mathcal{G}'_2 \subset \calo_{V_1}[W]$ is again finite, and the
construction of $\mbox{Coeff}\subset \calo_Z[W]$ is compatible
with base change. So $$\mbox{Coeff}(\mathcal{G}'_1)\subset
\mbox{Coeff}(\mathcal{G}'_2)\subset \calo_{Z_1}[W]$$ is the total
transform of $\mbox{Coeff}(\mathcal{G}_1)\subset
\mbox{Coeff}(\mathcal{G}_2)\subset \calo_{Z}[W]$.

By further restriction of $Z_1$ to $C$, we may assume that $Z_1$
is one dimensional. Theorem \ref{TH513} follows now from Prop
\ref{prop59}.
\endproof

\begin{theorem}\label{TH514} Let $\mathcal{G}(=G(\mathcal{G}))$ be a Diff-algebra over a smooth scheme $V$.
Fix a point $x\in \Sing(\mathcal{G})$, a smooth subscheme
$Z\subset V$ containing $x$, and two local (or formal)
retractions, say $\pi: V\to Z$ and $\pi': V\to Z$ at $x$. If
$\mbox{Coeff}(\mathcal{G})$ and $\mbox{Coeff}(\mathcal{G})'$ are
defined in terms of $\pi$ and $\pi'$ respectively, then both
define same Diff-algebra in $\calo_{Z,x}[W]$.
\end{theorem}
\proof Let $G(\mathcal{G})$ denote the Diff-algebra
spanned by $\mathcal{G}$ in the smooth scheme $V$. The claim is a
corollary of \ref{lem5885} .

\endproof

\section{ Further applications.}

There is a particular but natural morphism among smooth schemes, namely that defined by blowing up closed and smooth centers (i.e.,  monoidal transformations).  Given an ideal in a smooth scheme, there are several notions of transformations of sheaves of ideals defined in terms of monoidal transformations (e.g. total transforms,  weak transforms, and strict transforms of ideals.).

Questions as resolution of singularities, or Log principalization of ideals,  are formulated in terms of these notions of transformations. In the case of schemes over fields of characteristic zero, both resolution and Log principalization of ideals are two well known  theorems due to Hironaka.  If two ideals have the same integral closure, then a Log-principalization of one of them is also a Log-principalization of the other; the key point being that the transforms of both ideals also have the same integral closure.

These notions of transformations of ideals extend naturally to Rees algebras. And again, if two Rees algebras have the same integral closure, then  their transforms are Rees algebras with  the same
integral closure.

Both theorems of Log-principalization of ideals and of resolution of singularities are proved by induction on the dimension of the ambient space. In the setting of Diff-algebras this form of induction relates to the notion of restriction to a smooth subschemes, say $Z\subset V$ in Theorems \ref{TH513}.

The outcome of Theorems \ref{TH514} is that such form of restriction to $Z$ is, up to integral closure, independent of the particular retraction.  This result plays a role in the extension of resolution theorems to Rees algebras. 

Our development will be applied in \cite{VV2}, in relation with the study of hypersurface singularities  over fields of positive characteristic.


\smallskip


\begin{thebibliography}{W}

\bibitem{EncVil06} S. Encinas and O. Villamayor,  `Rees algebras and resolution of singularities'. To appear Rev. Mat.
Iberoamericana.

\bibitem{FN} M.~Fern\'andez-Lebr\'on and Luis Narv\'aez-Macarro, `Hesse-Schmidt Derivations and Coefficient Fields in Positive Characteristics',
 {\em J. Algebra}, 265 (2003), no. 1, 200--210.


\bibitem{Gir} J. Giraud, `Sur la
theorie du contact maximal',
     {\em Math. Zeit.}, 137  (1972), 285-310.




\bibitem{Giraud1975} J.~Giraud.
\newblock `Contact maximal en
caract\'eristique positive',
 \newblock {\em Ann. Scien. de l'Ec. Norm.
Sup. 4\`eme s\'erie}, 8
 (1975) 201-234.





\bibitem{Groth} A.~Grothendieck and J. Dieudonn\'e, `\'El\'ements de G\'eom\'etrie Alg\'ebrique IV:\'Etude locale des sch\'emas (Quatri\`eme Partie),
  vol. 32. Inst. Hautes \'Etudes Sci. Publ. Math., Press Univ. de France, Paris, (1967).




\bibitem{Hironaka77} H. Hironaka,
`Idealistic exponent of a
 singularity', {\em Algebraic Geometry},  The
John Hopkins centennial
 lectures, Baltimore, John Hopkins University Press
(1977), 52-125.

\bibitem{Hironaka03} H. Hironaka,
`Theory of infinitely near singular points', Journal Korean Math.
Soc. 40 (2003), No.5, pp. 901-920


\bibitem{Hironaka05} H. Hironaka,
`Three key theorems on infinitely near singularities ',
Singularit\'es Franco-Japonaises. Jean-Paul Brasselet - Tatsuo
Suwa. Seminaires et Congr\'es 10 SMF (2005)


\bibitem{kaw} H. Kawanoue, `Toward resolution of singularities over a field of positive characteristic' arXiv:math.AG/ 0607009 3 Aug 2006.


\bibitem{kollar} J. Koll\'{a}r, `Resolution of Singularities-
Seattle Lecture' arXiv:math.AG/ 0508332 v1  17 Aug 2005.


\bibitem{LT} J. Lipman and B. Teissier , `Pseudo-rational local rings and a theorem of Briancon-Skoda about integral closures of ideals'. Michigan Math. J.  (1981) 28, pp 97-116.




\bibitem{Matsumura} H. Matsumura, {\it Commutative Algebra (second Edition)}
Mathematics Lecture Note Series, Benjamin, Cumming Publishing
Company (1980).

\bibitem{MV} K. R. Mount and O. Villamayor (Sr) {\it Taylor Series and higher derivations.}
Vol 18, Impresiones del Departamento de Matematicas, Universidad
de Buenos Aires (1968, reimpresion 1979).

\bibitem{LN} L. Narv\'{a}ez Macarro {\em A note on the behaviour
under a ground field extension of quasicoefficient fields} J.
London Math. Soc. (2) 43 (1991) 12-22.

\bibitem{t} W. Traves `Localization of the Hasse-Schmidt algebra', Canad. Math. Bull. 46
(2003), no. 2, 304--309.



\bibitem{Vasconcelos} W. V. Vasconcelos {\it Arithmetic of Blowup
Algebras} London Mathematical Society. Lecture Notes Series 195.
Cambridge University Press 1994.

\bibitem{VVO}
O.~Villamayor, {\em Differential operators on smooth schemes and
embedded singularities.} Preprint.

\bibitem{VV2}
O.~Villamayor, {\em Differential operators and singularities in
positive characteristic. } In preparation.



 \bibitem{WL} J. Wlodarczyk, `Simple Hironaka resolution in characteristic
 zero', to appear (2005).

  \bibitem{Y} A. Yakutieli, `An explicit construction of the
  Grothendieck residual complex', Asterisque, vol. 208, Soc. Math.
  France, (1992).

   \bibitem{You} B. Youssin, `Newton Polyhedra without
   coordinates', Mem. Amer. Math. Soc. 433 (1990), 1-74,75-99.

 \bibitem{ZS} O. Zariski, P. Samuel {\it Commutative Algebra (Vol
 II)} University Series in Higher Mathematics. D. Van Nostrand
 Company (1960).


 \end{thebibliography}
\end{document}